\documentclass{article}
\usepackage{amsmath, amssymb,cmll,amsthm,nicefrac} %,amsthm,txfonts
\usepackage{multicol}
\usepackage{fourier}
\usepackage{bussproofs}
\usepackage{enumerate}
\usepackage{pgf}
\usepackage{tikz}
\usepackage{subcaption}
\usepackage{float}
\restylefloat{figure}
\newcommand{\bc}{\begin{C}}
\newcommand{\ec}{\end{C}}
\newcommand{\be}{\begin{equation}}
\newcommand{\ee}{\end{equation}}

\newcommand{\nb}{\begin{Prop}}
\newcommand{\nbe}{\end{Prop}}
\newcommand{\bl}{\begin{LE}}
\newcommand{\el}{\end{LE}}

\newcommand{\bd}{\begin{Def}}
\newcommand{\ed}{\end{Def}}
\newcommand{\bt}{\begin{Th}}
\newcommand{\et}{\end{Th}}

\newtheorem{Th}{Theorem}[section]
\newtheorem{LE}[Th]{Lemma}
\newtheorem{C}[Th]{Corollary}

\newcommand{\bp}{\begin{Prop}}
\newcommand{\ep}{\end{Prop}}

\newtheorem{Prop}[Th]{Proposition}

\theoremstyle{definition}
\newtheorem{Def}[Th]{Definition}

\theoremstyle{definition}
\newtheorem{Rem}[Th]{Remark}

\newcommand{\bbot}{\bot \!\!\!\!\!\! \bot}

\DeclareMathSymbol{\mlq}{\mathord}{operators}{``}
\DeclareMathSymbol{\mrq}{\mathord}{operators}{`'}

\begin{document}
 \title{Tensor term logic for categorial grammars: simple unification of commutative and noncommutative structure
}
\author{Sergey Slavnov
\\HSE University, Moscow, Russia
\\{slavnovserge@gmail.com}}
 \maketitle
\begin{abstract}
A prototypical example of  categorial  grammars are  those based on Lambek calculus, i.e. noncommutative intuitionistic  linear logic. However, it has been noted that purely noncommutative operations are often not sufficient for modeling even very simple natural language phenomena.
Therefore a number of alternative formalisms are considered in the literature: those using purely ``commutative'' linear logic as well as combining (to some level) commutative and non-commutative operations. The logic of tensor terms that we propose is a variant of such a combination. This logical calculus was designed specially for defining categorial grammars. It contains Lambek calculus and multiplicative linear logic as conservative fragments, yet the syntax is sufficiently simple, not departing much from that of multiplicative linear logic. The system is cut-free and decidable, it has both intuitionistic and classical versions, and is equipped with a simple intuitive semantics, which is sound and complete.
\end{abstract}

\section{Introduction}
A prototypical example of  categorial  grammars are  those based on Lambek calculus ({\bf LC}) \cite{Lambek}, \cite{Lambek_empty} or on some of its fragments. (Although, historically, first categorial grammars \cite{Ajdukiewicz}, \cite{BarHillel} were written earlier than {\bf LC} was formulated.)

{\bf LC} itself, from the mathematical point of view, is nothing but noncommutative intuitionistic  linear logic. However, it has been noted that purely noncommutative operations are often not sufficient for modeling even very simple natural language phenomena.
Therefore a number of alternative formalisms are considered in the literature: those using purely ``commutative'' linear logic as well as combining (to some level) commutative and non-commutative operations. The {\it logic of tensor terms} that we propose in this paper is a variant of such a combination. For motivation, let us discuss in some greater detail why commutativity in a categorial grammar might be wanted.

\subsection{When noncommutativity becomes a problem}\label{first label}
Noncommutativity of the underlying logic accounts for word order in the grammar and thus seems a natural, necessary property. However, in some situations in poses problems.
A simple example is formation of relative clause with medial extraction.

 Typically, for constructing an object relative clause in English using {\bf LC}, one adds to the lexicon a relative pronoun of type $(NP\backslash NP)/(S/NP)$ (where $NP$ stands for noun phrases, and $S$ for sentences.) Having the lexicon
 $$
 \mbox{John},\mbox{Mary}:NP,\quad \mbox{loves}:(NP\backslash S)/NP,\quad\mbox{whom}:(NP\backslash NP)/(S/NP),
 $$
 we derive that for any $x$ of type $NP$ the string ``John loves $x$'' is of type $S$, hence
 ``John loves'' is of type $S/NP$. This corresponds to the {\bf LC} derivation
 $$
 \cfrac{NP,(NP\backslash S)/NP,NP\vdash S}{NP,(NP\backslash S)/NP\vdash S/NP}(/\rm{R}).
 $$
  Then the concatenation of ``Mary'', ``whom'' and ``John loves'' must have type $NP$, i.e. the string ``Mary whom John loves'' is a noun phrase as desired.

  However it does not work that well when we add to the lexicon the adverb ``madly'' (of type $(NP\backslash S)\backslash(NP\backslash S)$) and want to derive that the string ``Mary whom John loves madly'' is of type $NP$ as well. We easily get that for any $x$ of type $NP$ the string
  ``John loves $x$ madly'' is of type $S$: this corresponds to the {\bf LC} derivable sequent
 \be\label{john loves x madly}
 {NP,(NP\backslash S)/NP,NP,(NP\backslash S)\backslash(NP\backslash S)\vdash S}.
  \ee
  But we cannot somehow extract this ``$x$'' from the sentence to create the relative clause. The ``$x$'' in question corresponds to the formula $NP$ to the left of the turnstile in (\ref{john loves x madly}), and that formula, being surrounded by other formulas on both sides, cannot be moved to the right of the turnstile precisely because the logic is noncommutative.
 Intuitively, we might want to say that ``John loves(...)madly'' is a  ``sentence with a hole'' and has type $NP\multimap S$, where $\multimap$ is the usual ``commutative'' implication of linear logic, and that ``whom''  accepts objects of this type as one of its arguments. But, since the usual commutative implication is missing in {\bf LC}, the noun phrase ``Mary whom John loves madly'' remains underivable.

 There is also another, very different class of situations, where commutativity might be wanted. Namely, it is the case of {\it scrambling} or the so called ``free word order''. We do not consider this problem in the current work.

 \subsection{Commutative extensions and variations}
  One of the approaches to restoring commutative operations is to consider {\it multimodal categorial grammars} (see \cite{Moortgat}). These use rather complicated extensions of {\bf LC} that have multiple {\it modes of composition}. That is, instead of a single family $\{\bullet,\backslash,/\}$ of {\bf LC} connectives we have multiple families
  $\{\bullet_i,\backslash_i,/_i\}$, where $i$ is the ``mode index''.
  Some of these modes might be commutative and some noncommutative (and also associative and nonassotiative), and they may interact with each other in nontrivial ways. Sequent antecedents in such a system are neither lists, nor multisets of formulas, but rather labeled binary trees, where nodes are labeled by modes of composition. Sometimes, systems of this sort may have also unary modal operators (and unary branches in the antecedent trees).

 An alternative is to change ``elementary building blocks'' used for constructing phrases.
  {\it Abstract categorial grammars} ({\bf ACG}) \cite{deGroote} (see also  \cite{Muskens}, \cite{PollardMichalicek}) construct phrases from {\it linear $\lambda$-terms}. (Indeed, $\lambda$-calculus is a programming language,  and the string data type together with all necessary operations can be represented in it). The natural typing system for linear $\lambda$-terms is (ordinary, commutative) multiplicative intuitionistic linear logic
 ({\bf MILL}), so the underlying logical calculus for {\bf ACG} is fully commutative. The word order is encoded in the ``low-level'' $\lambda$-terms. Unfortunately,  there are some linguistic situations (for example, non-constituent coordination) where, on the contrary, the noncommutative structure of {\bf LC} is wanted and such systems  as {\bf ACG} behave very poorly, see \cite{Moot_inadequacy} for a discussion. One could also remark that, although the underlying commutative logic is familiar and intuitive, reading concrete strings (``surface forms'') from abstract $\lambda$-terms is sometimes a complicated task.

 {\it Hybrid type logical categorial grammars} ({\bf HTLCG}) \cite{KubotaLevine} extend {\bf ACG} by adding to the syntax of $\lambda$-terms an explicit operation corresponding to string concatenation and having both {\bf MILL} and {\bf LC}  operations on types.
There are also other logical systems where commutative and noncommutative connectives coexist, such as
 {\it non-commutative logic}  \cite{AbrusciRuet} and  {\it partially commutative logic} \cite{Partially_commutative}. An apparent drawback of such systems is that combining two sorts of connectives leads to rather complicated syntax (and unclear semantics).

  Finally {\it displacement calculus} ({\bf DC}) \cite{Morrill_Displacement} is an extension of {\bf LC} where elementary building blocks are ordered string tuples and there are infinitely many connectives corresponding to different ways of putting such tuples together. Unfortunately, the system also suffers from the  complicated syntax. In fact, even a sequent calculus in the usual sense is absent, and for a cut-free formulation one has to use specific complicated structure of {\it hyper-sequents}.

\subsection{First order linear logic and categorial grammars}
 A remarkable observation \cite{MootPiazza} is that  first order (commutative) multiplicative linear logic ({\bf MLL1} or {\bf MILL1} for the intuitionistic version) can be used as a kind of grammatical formalism. In particular, both {\bf LC} and {\bf ACG}, as well as {\bf HTLCG}, can be encoded as conservative fragments of {\bf MILL1} \cite{Moot_extended}, \cite{Moot_comparing}, \cite{Moot_inadequacy}. (A similar encoding for {\bf DC} is also available, but in this case the embedding is not conservative, see \cite{Moot_comparing}.)
 Usage of {\bf MLL1} or {\bf MILL1} might be
 an alternative approach to combining commutative and noncommutative structure, avoiding syntactic complications and staying in the limits of familiar commutative linear logic and usual sequent calculus.

   The encodings in question are based on interpreting  first order variables as denoting {\it string positions}. Noncommutative constructions, such as string  concatenation, are represented by means of quantifiers binding variables. Typically, an {\bf LC} formula $F$ has the {\bf MILL1} translation $||F||^{x,y}$ parameterized by the free variables $x,y$ defined by induction as
   \be\label{Lambek2MILL1}
   \begin{array}{rl}
||p||^{x,y}=p(x,y),&|||A\bullet B||^{x,y}=\exists t(||A||^{x,t}\otimes||B||^{t,y}),\\
||A\backslash B||^{x,y}=\forall t(||A||^{t,x}\multimap||B||^{t,y}),&
||B/ A||^{x,y}=\forall t(||A||^{x,t}\multimap||B||^{y,t}),
\end{array}
   \ee
   where $p$ is a literal.

   However, one might argue that the full system of {\bf MILL1} is quite redundant for the above linguistic interpretations.
   Actually, only a tiny fragment  of the system is used in these constructions, and the ``linguistic'' meaning of general {\bf MILL1} formulas and sequents seems,
at least, obscure. Moreover,
   the fragment of interest is not closed either under subformulas or proof-search. Typically, when deriving in {\bf MILL1} the translation of an {\bf LC} sequent, in general, we have to use at intermediate steps sequents  that  have no clear linguistic meaning, if any.

  In \cite{Slavnov_generating} (elaborating on \cite{Slavnov_cowordisms}, \cite{Slavnov_tensor}) we proposed to isolate the relevant fragment of {\bf MILL1}, calling it the {\it linguistic fragment}, and then gave an alternative formulation for the latter. In a greater detail, we proposed
 the so called {\it extended tensor type calculus} ({\bf ETTC}), which is
  a certain calculus of typed terms, called {\it tensor terms}, and is
  equipped with an intrinsic deductive system and a transparent linguistic interpretation.
  Then we established strict  equivalence of the proposed calculus and the linguistic  fragment.
   %the

  In the current work, which can be seen as an  upgrade of the above, we propose the {\it logic} of tensor terms.  {\it Tensor term logic} ({\bf TTL}), is very similar to the previously defined {\bf ETTC}, but has much simpler syntax and sequent calculus, not too different from the ordinary propositional {\bf MLL}. Besides, {\bf TTL} is equipped with a sound and complete {\it semantics}, similar to the phase semantics of linear logic.
   (On the other hand,  the direct translation to the linguistic fragment of {\bf MILL1} in the case of {\bf TTL} is lost. Connections of {\bf TTL} to first order linear logic remain to be understood; we make several side remarks on this subject throughout the text.)

\subsection{Tensor term logic}
   {\it Tensor terms} \cite{Slavnov_tensor}, \cite{Slavnov_generating} are, basically,  tuples of words (written multiplicatively, as formal products) with labeled endpoints. We write words in square brackets and represent the endpoints as lower and upper indices,  lower indices standing for right endpoints and upper indices for left ones. Thus, tensor terms can have the form
$$[\mbox{a}]_i^j,\quad [\mbox{a}]_i^j\cdot[\mbox{b}]_k^l\cdot[\mbox{c}]_r^s,\quad [\epsilon]_i^j $$
(where $\epsilon$ stands for the empty word) and so on.

An index in a term can be repeated at most twice, once as an upper, and once as a lower one. A repeated  index is called {\it bound}; it means that the corresponding words are concatenated along matching endpoints. For example, we have the term equality $$[\mbox{a}]_j^i\cdot[\mbox{b}]_l^k\cdot[\mbox{c}]_k^j=[\mbox{acb}]_l^i$$
(the product is commutative). Naturally, indices that are not bound are called {\it free}.

In  the logic of tensor terms ({\bf TTL}), formulas or {\it types} denote sets of terms, and sequents denote relations between those (just as in {\bf LC} formulas denote sets of strings and sequents denote relations between strings). The syntax of {\bf TTL} is just that of {\bf MLL} with the addition that formulas are decorated with indices and the sequent calculus is supplemented with a couple of rules for index manipulation.  Typically, an atomic {\bf TTL} formula $p^{i_1\ldots i_m}_{j_1\ldots j_n}$ denotes a set of tensor terms all whose elements have
as their set of upper, respectively lower  free indices the collection
$i_1,\ldots,i_m$, respectively ${j_1,\ldots, j_n}$ . Tensor product of types corresponds to   product of terms (just as  product in {\bf LC} corresponds to concatenation of strings), and implication  is the adjoint operation.

{\bf TTL} formulas can also have bound indices; these may occur in  compound formulas. For example, the formula $a^i_j\otimes b^j_k$ denotes the product of term sets corresponding to $a^i_j$ and $b^j_k$. The latter sets have elements of the form $[u]^i_j$ and $[v]^j_k$ respectively; these are just strings with marked endpoints. Since  the endpoints match, the product becomes concatenation, and the denotation of $a^i_j\otimes b^j_k$
consists of elements of the form $[uv]^i_k$, with only two free indices. We say that the index $j$ in $a^i_j\otimes b^j_k$ is {\it bound}.

The above sketchy analysis suggests how to encode {\bf LC} to {\bf TTL}. In fact, the encoding is almost the same as the translation to {\bf MILL1} in (\ref{Lambek2MILL1}), with the only difference  that we do not need cumbersome quantifiers.
 The translation of an {\bf LC} formula $F$ is parameterized by two free indices according to the rules
  \[
   \begin{array}{rl}
||p||^i_j=p^i_j,&|||A\bullet B||^i_j=||A||^i_k\otimes||B||^k_j,\\
||A\backslash B||^i_j=||A||^k_i\multimap||B||^k_j,&
||B/ A||^i_j=||A||^j_k\multimap||B||^i_k.
\end{array}
\]

On the other hand, when the indices in subformulas are disjoint, we obtain discontinuous operations typical for ``commutative'' grammars such as {\bf ACG}. Tensor corresponds to ``disjoint union'' and implication produces a ``hole''in a string. In this way we get both ``commutative'' and ``noncommutative'' operations available, without departing much from the simple syntax of propositional {\bf MLL}. Cut elimination, decidability and a simple, intuitive semantics are also at hand.

  \section{Tensor terms}\label{tensor terms section}
  \bd\label{tensor terms}
Given  an alphabet $T$ of {\it terminal symbols} and an infinite set $\mathit{Ind}$ of {\it indices}, we consider the free commutative monoid
(written multiplicatively with the unit denoted as 1)
of {\it tensor pseudo-terms}  (over $T$ and $\mathit{Ind}$) generated by the sets
\be\label{generators}
\{[w]_i^j|~i,j\in\mathit{Ind},w\in T^*\}\mbox{ and }\{[w]|~w\in T^*\}.
\ee
 A pseudo-term $t$ is {\it well-formed} if   any index in $t$ has at most one lower and one upper occurrence.
 The set  of {\it tensor terms} (over $T$ and $\mathit{Ind}$) is the set of well-formed pseudo-terms
 quotiented
by the
 {\it congruence} generated by the relations
\be\label{tensor term relations}
\begin{array}{c}
[u]_j^i\cdot[v]_k^j\equiv[uv]_k^i,\quad [w]^i_i\equiv[w],
\\
 {[} a_1\ldots a_n ]
 \equiv[a_na_1\ldots a_{n-1}]\mbox{ for }a_1,\ldots,a_n\in T.
\end{array}
\ee

Generating set (\ref{generators}) elements of the form $[w]_i^j$ are {\it elementary regular tensor terms} and those of the form $[w]$ are elementary {\it singular tensor terms}.
\ed

 An index  occurring in a well-formed pseudo-term $t$ is {\it free} in $t$ if it occurs in $t$ exactly once. An index occurring in $t$ twice (once as an upper one and once as a lower one) is {\it bound} in $t$.
    Sets of free indices are invariant under congruence (\ref{tensor term relations}), so that  free indices of a tensor term are well-defined as well.
    \bd
    For a tensor term or a well-formed tensor pseudo-term $t$ we denote    the set of its free upper, respectively lower indices as
     $FI^\bullet(t)$, respectively $FI_\bullet(t)$. We call the  pair  $FI(t)=(FI^\bullet(t),FI_\bullet(t))$ the {\it boundary} of $t$.
%
%     A well-formed pseudo-term is {\it normal} if it has no bound indices.  A {\it normal form} of a tensor term $t$ is a normal pseudo-term representing $t$. When a term $t$ is written in a normal  form,
%  $t=[w_1]^{i_1}_{j_1}\ldots[w_k]^{i_k}_{j_k}[u_1]\ldots[u_l]$, where $i_1,\ldots, i_k$, $j_1,\ldots, j_k$ are pairwise distinct, the elementary terms
%  $[w_1]^{i_1}_{j_1},\ldots,[w_k]^{i_k}_{j_k},[u_1],\ldots,[u_l]$
%  are the {\it components} of $t$.
\ed

  The meaning of relations (\ref{tensor term relations}) is as follows. Terms stand for string tuples, indices denote string endpoints. Repeated, i.e. bound, indices correspond to matching endpoints and strings are glued along them.  The  unit $1$ corresponds to the empty tuple. As for singular terms, such as $[w]^i_i$,  they correspond to {\it cyclic} words (this explains the last relation in (\ref{tensor term relations})). These arise when strings are glued cyclically.
    Singular terms are pathological, but we need them for consistency of definitions.

  One can take a more geometric view and represent elementary regular terms as  directed  string-labeled edges with indices denoting vertices and
  elementary singular terms as closed loops labeled with cyclic words.
   Edges are glued along matching vertices. General terms become  graphs constructed from elementary ones by taking disjoint unions and
   gluing edges along matching vertices.  The title ``boundary'' is inspired precisely by such a view.
   %The geometric representation
%   makes it especially evident that
%    any tensor term   has unique normal form. The  components of a tensor term correspond to geometric components (edges and loops) comprising the corresponding graph.
\bd
%A tensor term is {\it regular} if all its components are regular, otherwise the term is {\it singular}.
For a well-formed pseudo-term $t$  with $i\in FI^\bullet(t)$ and $j$ not occurring in $t$, the {\it reparameterization}
$t^{[\nicefrac{j}{i}]}$ of $t$
is the  pseudo-term obtained from $t$ by replacing the index $i$ with $j$.    Similarly, for $t$ with $i\in FI_\bullet(t)$ and $j$ not occurring in $t$, the {\it reparameterization} $t_{[\nicefrac{j}{i}]}$ is $t$ where $i$ is replaced with $j$. The same convention applies  when $t$ is a tensor term.
\ed
It is easy to see that reparameterization  takes well-formed pseudo-terms to well-formed pseudo-terms and preserves congruence classes, hence it is indeed well-defined for terms as well.

We have the following {\it zoology} of tensor terms.
\bd
A tensor term $t$ is {\it regular} if it can be written as $t=[w_1]^{i_1}_{j_1}\cdots[w_n]^{i_n}_{j_n}$, where the indices $i_1,j_1,\ldots,i_n,j_n$ are pairwise distinct. The term $t$ is {\it singular} if it has the form $[u_1]\cdots[u_k]t'$, where $k>0$ and $t'$ is regular.

The term $t$ is {\it lexical} if it has the form $t=[u_1]\cdots[u_k]\cdots[w_1]^{i_1}_{j_1}\cdots[w_n]^{i_n}_{j_n}$, where
the indices $i_1,j_1,\ldots,i_n,j_n$ are pairwise distinct and all words $u_1,\ldots, u_k,w_1,\ldots,w_n$ are nonempty.

  The elementary term  $[\epsilon]^j_i$ where $i\not=j$ and $\epsilon$ denotes the empty word, is called {\it Kronecker delta} and  is denoted  as $\delta^j_i$.
 % A regular term is {\it lexical} if none of its components is a Kronecker delta. A lexical term is {\it basic lexical} if each its component has the form $[a]^i_j$, where $i$ is a single terminal symbol.
\ed
%Geometrically,  Kronecker deltas are edges carrying no labels.
% regular terms are edge-labeled bipartite graphs, lexical terms are bipartite graphs, where every edge carries a nonempty label, and basic lexical terms are bipartite graphs, where every edge is labeled with a terminal symbol.
\bp\label{zoology}
If $t$ is a  term  with
$i\in FI^\bullet(t)$, $j\in FI_\bullet(t)$  and
$i',j'$ are fresh then
 $\delta_{i}^{i'}\cdot t= t^{[\nicefrac{i'}{i}]}$ and
$\delta_{j'}^{j}\cdot t= t_{[\nicefrac{j'}{j}]}$.
 $\Box$
\ep
\begin{Rem}\label{Rem1}
In {\it extended tensor type calculus} ({\bf ETTC}) of \cite{Slavnov_generating} we had the additional term congruence
\be\label{empty loop}
[\epsilon]\equiv 1.
\ee
 This was needed for the main result of that paper, namely, for conservative translation of {\bf ETTC} to first order linear logic {\bf MLL1}. In the current work we drop this congruence, because it leads to unnecessary complications.
\end{Rem}
It should be  emphasized that multiplication of terms is only partially defined: the product $ts$ exists as a well-defined term only if $FI^\bullet(t)\cap FI^\bullet(s)=FI_\bullet(t)\cap FI_\bullet(s)=\emptyset$. It is associative in the following sense: if both $(ts)k$ and $t(sk)$ are well-defined then they are equal. However, it may be the case that one side is undefined, while the other one is  a well-defined term.

We will assume the following convention.
\bd\label{iterated multiplication}
The expression $t_{(1)}\ldots t_{(n)}$ denotes the iterated  product of terms satisfying the condition
$FI^\bullet(t_{(i)})\cap FI^\bullet(t_{(j)})=FI_\bullet(t_{(i)})\cap FI_\bullet(t_{(j)})=\emptyset$ for all $i\not=j$.
\ed
(In the definition above we put sequence numbers in brackets in order to avoid confusion with indices.)

 In the sequel, when writing  products of terms we will  implicitly assume that the operations are well-defined and,  in the case of iterated products, the convention from the above definition  applies. This  makes expressions with iterated products unambiguous.

\section{Intuitionistic tensor term logic}
\subsection{Intuitionistic tensor types}\label{int_tensor_types}
\bd
Given a set $\mathit{Lit}$ of   {\it literals}, where every  element $p\in \mathit{Lit}$ is assigned a {\it valency} $v(p)\in \mathbb{N}^2$ and an infinite set $\mathit{Ind}$ of indices,
the set $\widetilde{\mathit{Fm}}$ of {\it intuitionistic tensor pseudoformulas} over the alphabet $\mathit{Lit}$ is built  according to the grammar in Figure \ref{ITTP_lang}, where $\mathit{Prop}$ is the set of {\it propositional symbols}. The set $\widetilde{\mathit{At}}$ is the set of {\it atomic pseudoformulas}.
\ed
 The connectives $\otimes$ and $\multimap$ are called, respectively, {\it tensor} (also {\it times}) and {\it linear implication}. The propositional constant ${\bf 1}$ is called {\it one} (also the {\it unit}).

\bd
For every pseudoformula $A$ (over the given alphabet $\mathit{Lit}$) the set of index occurrences in $A$ is partitioned into the pair $(I^\bullet(A),I_\bullet(A))$ of sets, the set of  {\it upper} and {\it lower polarity} indices respectively.   The partition is defined by induction on $A$ in Figure \ref{ITTP_polarities}.

An atomic pseudoformula  $p^{i_1\ldots i_m}_{j_1\ldots j_n}$ is {\it well-formed}  if all indices $i_1,\ldots,i_m$, $j_1,\ldots,j_m$ are pairwise distinct. A {general} intuitionistic pseudoformula $A$ is {\it well-formed} if it is built of well-formed atomic ones and any index  has at most one  upper and one  lower polarity occurrence in $A$.  A well-formed intuitionistic tensor pseudoformula $F$ is an {\it intuitionistic tensor formula} or an {\it intuitionistic tensor type} (over $\mathit{Lit}$). A well-formed atomic intuitionistic tensor pseudoformula is an {\it atomic intuitionistic tensor formula}. The set of atomic formulas is denoted as $At$.
The set of intuitionistic tensor formulas over the given $\mathit{Lit}$ is the {\it intuitionistic tensor language generated by} $\mathit{Lit}$ or, simply the {\it intuitionistic tensor language over} $\mathit{Lit}$.
\ed
We emphasize that there is a certain notational inconsistency in our language. The polarity of an index occurrence does not necessarily coincide with its actual graphical location in the formula. Typically, upper indices written in the antecedent of an implication will have lower polarity and so on. Such an inconsistency does not arise when tensor formulas are written in the format of classical, rather than intuitionistic, linear logic.
\bd
 An index in a tensor formula is {\it free}  if it occurs exactly once, otherwise it is {\it bound}. The set of free upper, respectively lower indices of a tensor formula $A$ is denoted as $FI^\bullet(A)$, respectively $FI_\bullet(A)$. The pair $FI(A)=(FI^\bullet(A),FI_\bullet(A))$ is the {\it boundary } of $A$. Tensor formulas $F,G$  are $\alpha$-equivalent ($F\equiv_\alpha G$) if they can be obtained from each other by successively renaming  bound indices.
 \begin{figure}
  \centering
  \subfloat[Intuitionistic tensor pseudoformulas
  \label{ITTP_lang}
]
{
$\begin{array}{c}
  \mathit{Prop}=\mathit{Lit}\cup\{{1}\},
  \quad v({1})=(1,1),
 \\[.1cm]
%\quad
\widetilde{\mathit{At}}=\{p_{j_1\ldots j_n}^{i_1\ldots i_m}|~p\in \mathit{Prop},i_1\ldots i_m,j_1\ldots j_n\subset\mathit{Ind},v(p)=(m,n)\},
\\[.1cm]
\widetilde{\mathit{Fm}}::=\widetilde{\mathit{At}}|(\widetilde{\mathit{Fm}}\otimes\widetilde{\mathit{Fm}})|(\widetilde{\mathit{Fm}}\multimap\widetilde{\mathit{Fm}}).
\end{array}
$
}\\[.1cm]
 %%%%%%%%%%%%%%%%%%%%%%%%%%%%%%%%%%%%%%%%%%%
 \subfloat[Index polarities
  \label{ITTP_polarities}
]
{
$
\begin{array}{rl}
I^\bullet({p^{i_1,\ldots,i_n}_{j_1,\ldots,j_m}})=\{i_1,\ldots,i_n\},&
I_\bullet({p^{i_1,\ldots,i_n}_{j_1,\ldots,j_m}})=\{j_1,\ldots,j_m\}
\mbox{ for }p\in\mathit{Prop},\\
I^\bullet(A\otimes B)=I^\bullet(A)\cup I^\bullet(B),&
I_\bullet(A\otimes B)=I_\bullet(A)\cup I_\bullet(B),
\\
I^\bullet(A\multimap B)=I^\bullet(B)\cup I_\bullet(A),&
I_\bullet(A\multimap B)=I_\bullet(B)\cup I^\bullet(A).
\end{array}
$
}
\caption{{\bf ITTL} language}
\label{ITTP_language}
\end{figure}

 The notation $A^{[\nicefrac{i}{j}]}$ stands for $A$ with $i$ substituted instead of a free upper occurrence of $j$ (if there is any), where it is implicitly assumed that the result of the substitution is a well-formed tensor formula and the index  $i$ does not become bound after substitution, i.e. if $j\in FI^\bullet(A)$ then $i$ does not occur in $A$. Similarly,
    $A_{[\nicefrac{i}{j}]}$ stands for $A$ with $i$ substituted instead of a free lower occurrence of $j$, with the implicit assumption that if $j\in FI_\bullet(A)$ then $i$ does not occur in $A$. When two formulas differ by such a substitution of free indices we say that they are {\it reparameterizations} of each other. Two formulas that can be obtained from each other by a chain of reparameterizations and $\alpha$-equivalences are said to have the same {\it  type symbol}. The {\it valency} of a type symbol $a$ is the pair $(\#FI^\bullet(A),\#FI_\bullet(A))$, where $A$ is any formula with symbol $a$ and $\#$ denotes the cardinality of a finite set.  A type symbol is {\it balanced} if its valency is of the form $(n,n)$. An intuitionistic tensor language  is balanced if all its literals are balanced.
\ed
The intended semantics is that tensor types are sets of tensor terms. In a greater detail, each type $A$ should be assigned a set $||A||$ of terms with boundary $FI(A)$, and this assignment should be consistent with reparameterization: for $i\in FI^\bullet(A)$ we have $t\in ||A||$ iff $t^{[\nicefrac{i}{j}]}\in||A^{[\nicefrac{i}{j}]}||$, similarly with lower indices.

A ``naive'' semantics of connectives is as follows. The tensor product $A\otimes B$ stands for the set $\{ts|~t\in ||A||,s\in ||B||\}$, and for the implicational type $A\multimap B$ we have $t\in||A\multimap B||$ iff for all $s\in ||A||$ it holds that $ts\in ||B||$. As for the unit we may put $||{\bf 1}^i_j||=\{\delta^i_j\}$. Note that consistency with reparameterization implies that if $j\in FI^\bullet(A)$ and $i$ does not occur in $A$ then
\be\label{consistency with reparameterization}
||A^{[\nicefrac{i}{j}]}||=||{\bf 1}^{i}_j\otimes A||=||{\bf 1}^j_i\multimap A||,
\ee
 similarly for lower indices (this can be seen from Proposition \ref{zoology}). Such an interpretation is sound for the system of {\it intuitionistic tensor term logic} ({\bf ITTL}) that we are about to introduce, although its completeness is under question. In any case, it does not make much sense when the language is not balanced, because for a ``non-balanced'' formula $A=a^{i_1\ldots i_m}_{j_1\ldots j_n}$, where $m\not=n$ we cannot have any terms with boundary $FI(A)$.
A complete semantics for {\bf ITTL} in a balanced language, which is only a slight refinement of the above, will be discussed later.

\subsection{Sequent calculus}
\begin{figure}
\centering
%%%%%%%%%%%%%%%%%%%%%%%%
\subfloat[Defining rules
  \label{ITTL_rules}
]
{
$
\begin{array}{c}
\begin{array}{rl}
\quad\cfrac{A\in At}{A\vdash A}{\rm{(Id)}}
&
\cfrac{\Gamma\vdash A\quad A,\Delta\vdash B}{\Gamma,\Delta\vdash B}({\rm{Cut}})
\\
\cfrac{\Gamma\vdash A\quad B,\Delta\vdash C}{\Gamma,A\multimap B,\Delta\vdash C}({\rm{L}}\multimap)
&
\cfrac{\Gamma,A\vdash B}{\Gamma\vdash A\multimap B}({\rm{R}}\multimap)
\\
\cfrac{\Gamma, A,B\vdash C}{\Gamma,A\otimes B\vdash C}({\rm{L}}\otimes)
&
\cfrac{\Gamma\vdash A\quad\Delta\vdash B}{\Gamma,\Delta\vdash A\otimes B}({\rm{R}}\otimes)
\\
\cfrac{\Gamma,A\vdash B~ j\in FI_\bullet(A)}{\Gamma,A_{[\nicefrac{i}{j}]},{\bf 1}_j^i\vdash B}({\rm{L}}{1}_\to)
&
\cfrac{\Gamma\vdash A~ j\in FI^\bullet(A)}{{\bf 1}^i_j,\Gamma\vdash A^{[\nicefrac{i}{j}]}}({\rm{R}}{1}^\to)
\\
\cfrac{\Gamma,{\bf 1}^i_j,A\vdash B~ j\in FI^\bullet(A)}{\Gamma,A^{[\nicefrac{i}{j}]}\vdash B}({\rm{L}}{1}^\leftarrow)
&
\cfrac{\Gamma,{\bf 1}^i_j\vdash A~ j\in FI_\bullet(A)}{\Gamma\vdash A_{[\nicefrac{i}{j}]}}({\rm{R}}{1}_\leftarrow)
\end{array}\\[.6cm]
\mbox{In }({\rm{L}}{1}_\to),({\rm{R}}{1}^\to),({\rm{L}}{1}^\leftarrow),({\rm{R}}{1}_\leftarrow)~
i\not\in FI^\bullet(A)\cup FI_\bullet(A)
\end{array}
$
}
\\
%%%%%%%%%%%%%%%%%%%%%%%%%%%%%%%%%
\subfloat[Admissible rules
  \label{ITTL_admissible rules}
]
{
$
\begin{array}{c}
\begin{array}{rl}
\cfrac{\Gamma,A\vdash B~ j\in FI^\bullet(A)}{\Gamma,{\bf 1}^j_i,A^{[\nicefrac{i}{j}]}\vdash B}({\rm{L}}{1}^\to)
&
\cfrac{\Gamma\vdash A~ j\in FI_\bullet(A)}{\Gamma,{\bf 1}_i^j\vdash A_{[\nicefrac{i}{j}]}}({\rm{R}}{1}_\to)
\\
\cfrac{\Gamma,A,{\bf 1}_i^j\vdash B~ j\in FI_\bullet(A)}{\Gamma,A_{[\nicefrac{i}{j}]}\vdash B}({\rm{L}}{1}_\leftarrow)
&
\cfrac{{\bf 1}_i^j,\Gamma\vdash A~ j\in FI^\bullet(A)}{\Gamma\vdash A^{[\nicefrac{i}{j}]}}({\rm{R}}{1}^\leftarrow)
\end{array}\\[.6cm]
\mbox{In all rules above }
i\not\in FI^\bullet(A)\cup FI_\bullet(A)
\end{array}
$
}
\caption{{\bf ITTL} sequent rules}
\end{figure}
\bd
An {\it intuitionistic tensor context} is a finite set $\Gamma$ of intuitionistic tensor formulas such that for any distinct occurrences $A,B\in\Gamma$ it holds that
\be\label{indices in context}
 FI_\bullet(A)\cap FI_\bullet(B)=FI^\bullet(A)\cap FI^\bullet(B)=\emptyset.
 \ee
  We denote $FI^\bullet(\Gamma)=\bigcup\limits_{A\in\Gamma} FI^\bullet(A)$,
$FI_\bullet(\Gamma)=\bigcup\limits_{A\in\Gamma} FI_\bullet(A)$.

An {\it intuitionistic tensor sequent} is an expression of the form $\Gamma\vdash A$, where $\Gamma$ is an  intuitionistic tensor context, $A$  an intuitionistic tensor formula and it holds that
$$
FI^\bullet(\Gamma)\cap FI_\bullet(A)=FI_\bullet(\Gamma)\cap FI^\bullet(A)=\emptyset,\quad
FI^\bullet(\Gamma)\cup FI_\bullet(A)=FI_\bullet(\Gamma)\cup FI^\bullet(A).
$$
The {\it polarity of a free index occurrence $i$ in the sequent} $\Gamma\vdash A$ is the same as the polarity of $i$ in $A$ if $i$ is located in $A$ and is the opposite of its polarity in $\Gamma$ if it is located in $\Gamma$.
\ed
The definition says that every index occurring in a sequent  has {\it exactly}  one upper and  one lower polarity free occurrence (keeping in mind that polarities to the left of the turnstile are read upside down.)

The intended semantics of the sequent $A_{(1)},\ldots,A_{(n)}\vdash B$ (sequence numbers in brackets in order to avoid confusion with indices) is the following: if we have terms $t_{(1)},\ldots,t_{(n)}$ of types $A_{(1)},\ldots,A_{(n)}$ respectively, then the product $t_{(1)}\cdots t_{(n)}$ is of type $B$. (When $n=0$, i.e. the sequent has the form $\vdash B$, this should be understood that the term $1$ is of type $B$.)
\bd
The system of {\it intuitionistic tensor  term logic} ({\bf ITTL}) is given by the rules in Figure \ref{ITTL_rules}, where  it is assumed  that all expressions are well-formed sequents,
  i.e. each free index has exactly one upper and  one lower polarity  occurrence .
  \ed
  We encourage the reader to check that the ``naive'' semantics of types discussed in the preceding subsection is sound for the above system. The last four rules (i.e. the part that is not the usual intuitionistic linear logic) express the requirement that the interpretation of types should be consistent with reparameterization as in (\ref{consistency with reparameterization}). In this group, the rules
   marked with the arrow directed from the unit are {\it expansion} rules, and those with the arrow directed towards the unit are {\it reduction} rules.
   There are other, admissible,  expansion and reduction rules of a similar form and with the same meaning shown in Figure \ref{ITTL_admissible rules}. Their admissibility will be proven in the next subsection. Finally, note also that there is another rule introducing the unit, namely the special instance ${\bf 1}^i_j\vdash {\bf 1}^i_j$ of the $({\rm{Id}})$ axiom.
   \begin{Rem}\label{unit_remark}
     We call the propositional constant {\bf 1}  the {\it unit}, but this might be somewhat misleading.

     Indeed,  there exists a kind of monoidal unit in the algebra of tensor terms, namely the term $1$ that corresponds to the empty product of elementary terms and satisfies the equality $1\cdot t=t$ for all terms $t$. However the formula ${\bf 1}^i_j$ stands for  the Kronecker delta $\delta^i_j$ rather than for 1. In principle, one  might add to the syntax of the logic the {\it monoidal unit}  formula $\mathbb{1}$ (with no indices)
      standing for the term 1. The rules for the monoidal unit then will be the usual rules of intuitionistic linear logic with unit:
     $\cfrac{\Gamma\vdash A}{\mathbb{1},\Gamma\vdash A}$ and ${\vdash\mathbb{1}}$.

     Our unit {\bf 1} and the corresponding formulas ${\bf 1}^i_j$ should not be confused with the usual monoidal unit $\mathbb{1}$ of linear logic. In particular, the monoidal unit $\mathbb{1}$ should have valency $v(\mathbb{1})=(0,0)$, and for our unit ${\bf 1}$ we have $v({\bf 1})=(1,1)$. (In some sense,  the unit  of {\bf ITTL} is an incarnation of the  unit of Lambek calculus.)
   \end{Rem}
   \begin{Rem}\label{Rem2}
 It might be interesting to study what will be the logic of tensor terms when we add to the definition of terms the additional congruence
  from (\ref{empty loop}) (see Remark \ref{Rem1}). For example the ``naive'' semantics will validate the
  principle $\vdash (A^i_j\multimap A^i_k)\otimes(A^l_k\multimap A^l_j)$.
  (Indeed, we have $\delta^j_k\in ||A^i_j\multimap A^i_k||$, $\delta_j^k\in ||A^l_k\multimap A^l_j||$, and $\delta^j_k\delta^k_j$ becomes $1$ under (\ref{empty loop}).)  But the latter sequent is easily seen to be underivable in the given formulation of {\bf ITTL} without the Cut rule. (The sequent has no free indices, therefore the last rule in its derivation cannot be $({\rm{R}}{1}_\leftarrow)$. So it could only be the $(\otimes)$ rule, and this does not work.) That the proposed formulation of {\bf ITTL} is cut-free will be shown shortly.
  \end{Rem}

\subsection{Simplest properties}
\bd
Let $\Gamma\vdash A$ be a sequent with a free index $i$. A {\it similarity transformation} of $\Gamma\vdash A$ consists in replacing both free occurrences of $i$ with a fresh index. Two sequents are {\it similar} (notation $\sim$) if they can be obtained from each other by a chain of similarity transformations. {\it $\alpha$-Equivalence} of tensor sequents (notation $\equiv_\alpha$) is the smallest equivalence relation such that similar sequents and sequents differing by renaming bound indices in formulas are $\alpha$-equivalent.
\ed
\bp\label{ITTL properties}
In {\bf ITTL}
 \begin{enumerate}[(i)]
% \item In a derivable sequent, any free index has exactly one upper and one lower polarity free occurrence.
 \item Similar sequents are cut-free derivable from each other;
 \item rules in Figure \ref{ITTL_admissible rules} are admissible and can be emulated without use of the Cut rule;
    %\item for any formula $F$ the sequent $\vdash\overline{F},F$ is derivable;
      \item $\alpha$-equivalent formulas are provably equivalent, i.e. if $F\equiv_\alpha F'$ then the sequents $F\vdash F'$, $F'\vdash F$ are derivable;
      \item $\alpha$-equivalent sequents are derivable from each other;
  \item
 if   $j\in FI^\bullet(A)$, and $i$ is an index not occurring in $A$ then
  there are provable equivalences
  ${\bf 1}^i_j\otimes A \cong {\bf 1}_i^j\multimap A\cong A^{[\nicefrac{i}{j}]}$.

  Similarly if   $j\in FI_\bullet(A)$, and $i$ is an index not occurring in $A$ then
  there are provable equivalences
  $A\otimes {\bf 1}^j_i\cong {\bf 1}^i_j\multimap A\cong A_{[\nicefrac{i}{j}]}$.
 \end{enumerate}
\ep
{\bf Proof}
 \begin{enumerate}[(i)]
\item Assume that the sequent $\Gamma'\vdash A'$ is obtained from $\Gamma\vdash A$  by a similarity transformation replacing  index $j$ with $i$.  Applying to $\Gamma\vdash A$ an expansion rule we replace the upper polarity occurrence of $j$  with $i$ and introduce into the sequent the formula ${\bf 1}^i_j$. Then we use a reduction rule and replace the lower polarity occurrence of $j$ and get the primed sequent.
\item As an example consider the $({\rm{L}}{1}_\leftarrow)$ rule. Assume that we have a sequent of the form
$\Gamma,A,{\bf 1}_i^j\vdash B$, where $j\in FI_\bullet(A)$, $i\not\in FI^\bullet(A)$. By definition of intuitionistic tensor sequents there must be also an upper polarity free occurrence of $i$ located in some formula $F$, which is necessarily distinct from $A$.
%(otherwise the formula occurrence of $i$ in $A$ would be bound).
Assume for definiteness that $F$ is in $\Gamma$, so that $\Gamma=\Gamma',F$. Note that $j\not\in FI_\bullet(F)$ because we already have the lower polarity occurrence of $j$ in $A$. The derivation is shown in Figure \ref{L1<- rule}. Another case is $F=B$, which is treated similarly.
\item Induction on $F,F'$, the base cases being $F=A\otimes B$, $F'=A^{[\nicefrac{i}{j}]}\otimes B_{[\nicefrac{i}{j}]}$ or
$F'=A_{[\nicefrac{i}{j}]}\otimes B^{[\nicefrac{i}{j}]}$ and analogously when the main connective is implication. The derivation for the first case is shown in Figure \ref{alpha-equiv_tensor}. Other cases are treated analogously.
\item Follows from the two preceding clauses.
\item As an example we derive the equivalence ${\bf 1}^i_j\otimes A\cong A^{[\nicefrac{i}{j}]}$ for $j\in FI^\bullet(A)$, $i$ fresh in Figure \ref{reparameterization derivation}. $\Box$
\end{enumerate}

The crucial property of cut-elimination will be proven separately in Section \ref{cut-elimination section}, after the introduction of {\it classical} tensor term logic.
  \begin{figure}%[htb]
\centering
%%%%%%%%%%%%%%%%%%%%%%%%%%%%%
\subfloat[Emulating the
$({\rm{L}}{1}_\leftarrow)$ rule
\label{L1<- rule}
]
{
$
\def\fCenter{\ \vdash\ }
\AxiomC{$\Gamma',F,A,{\bf 1}_i^j\fCenter B~j\in FI_\bullet(A), i\in FI^\bullet(F)$}
\RightLabel{(${\rm{L}}{1}^\leftarrow$)}
\UnaryInf$\Gamma',F^{[\nicefrac{j}{i}]},A\fCenter B$
\RightLabel{(${\rm{L}}{1}_\to$)}
\UnaryInf$\Gamma',F^{[\nicefrac{j}{i}]},A_{[\nicefrac{i}{j}]},{\bf 1}^i_j \fCenter B $
\RightLabel{(${\rm{L}}{1}^\leftarrow$)}
\UnaryInf$\Gamma',F,A_{[\nicefrac{i}{j}]}\fCenter B$
\DisplayProof
$
}
%\\[1cm]
%%%%%%%%%%%%%%%%%%%%%%%%%%%%%%%%%%%%%%%%%%%%%
\subfloat[
Renaming bound indices
\label{alpha-equiv_tensor}
]
{
$
\def\fCenter{\ \vdash\ }
\AxiomC{$A\fCenter A$}
\AxiomC{$B\fCenter B$}
\RightLabel{(${\rm{R}}\otimes$)}
\BinaryInf$A,B\fCenter A\otimes B$
\RightLabel{(${\rm{L}}{1}_\to$)}
\UnaryInf$A, B_{[\nicefrac{i}{j}]},{\bf 1}^i_j\fCenter A\otimes B$
\RightLabel{(${\rm{L}}{1}^\leftarrow$)}
\UnaryInf$ A^{[\nicefrac{i}{j}]}, B_{[\nicefrac{i}{j}]}\fCenter A\otimes B$
\RightLabel{(${\rm{L}}\otimes$)}
\UnaryInf$ A^{[\nicefrac{i}{j}]}\otimes B_{[\nicefrac{i}{j}]}\fCenter A\otimes B$
\DisplayProof
$
}
\\[1cm]
%%%%%%%%%%%%%%%%%%%%%%%%%%%%%%%%%%%%%%%%%%%%%
\subfloat[
Reparameterization
\label{reparameterization derivation}
]
{
$
\def\fCenter{\ \vdash\ }
\AxiomC{${\bf 1}^i_j\fCenter {\bf 1}^i_j$}
\AxiomC{$A\fCenter A~j\in FI^\bullet(A)$}
\RightLabel{(${\rm{R}}\otimes$)}
\BinaryInf${\bf 1}^i_j,A,\fCenter{\bf 1}^i_j\otimes A$
\RightLabel{(${\rm{L}}{1}^\leftarrow$)}
\UnaryInf$A^{[\nicefrac{i}{j}]}\fCenter {\bf 1}^i_j\otimes A$
\DisplayProof
$
$\quad$
$
\def\fCenter{\ \vdash\ }
\AxiomC{$A\fCenter A~j\in FI^\bullet(A)$}
\RightLabel{(${\rm{R}}{1}^\to$)}
\UnaryInf${\bf 1}^i_j,A\fCenter A^{[\nicefrac{i}{j}]}$
\RightLabel{(${\rm{L}}\otimes$)}
\UnaryInf$ {\bf 1}^i_j\otimes A\fCenter A^{[\nicefrac{i}{j}]}$
\DisplayProof
$
}
\caption{Some {\bf ITTL} derivations}
\end{figure}

\section{Classical format}
\subsection{The logic}
\bd
Given a set $\mathit{Lit}_+$ of   {\it positive literals}, where every  element $p\in \mathit{Lit}_+$ is assigned a {\it valency} $v(p)\in \mathbb{N}^2$,
and an infinite set $\mathit{Ind}$ of indices,
the set $\widetilde{\mathit{Fm}}$ of {\it classical tensor pseudoformulas} over $\mathit{Lit}_+$  is built  according to the grammar in Figure \ref{TTL_lang},
where  $\mathit{Lit}_-$ and  $\mathit{Lit}$ are, respectively, the set of {\it negative literals} and of all literals and $\mathit{Prop}$ is the set of {\it propositional symbols}. The convention for negative literals is that $v(\overline{p})=(m,n)$ if  $v(p)=(n,m)$.
The set $\widetilde{\mathit{At}}$ is the set of classical atomic pseudoformulas.
\ed
  The connective $\parr$  is called   {\it cotensor} (also {\it par}). The constant  $\bot$ is the {\it counit}.
Duality $\overline{(.)}$ is not a connective or operator, but is definable.
Unlike the intuitionistic case we do not have to define  upper/lower index polarities different from actual index locations.
 \bd
An atomic tensor pseudoformula  $p^{i_1\ldots i_m}_{j_1\ldots j_n}$
(over the given alphabet $\mathit{Lit}_+$)
is {\it well-formed}  if all indices $i_1,\ldots,i_m$, $j_1,\ldots,j_m$ are pairwise distinct. A {general} classical tensor pseudoformula $A$ is {\it well-formed} if it is built from well-formed atomic ones and any index  has at most one  upper and one  lower occurrence in $A$.  A well-formed classical tensor pseudoformula $F$ is a {\it classical tensor formula} or a {\it classical tensor type}
(over  $\mathit{Lit}_+$). A well-formed classical atomic tensor pseudoformula is a {\it classical atomic tensor formula}. The set of classical atomic tensor formulas is denoted as $At$.
The set of classical tensor formulas over the given alphabet $\mathit{Lit}_+$ is the {\it classical tensor language generated by} $\mathit{Lit}_+$ or, simply the {\it classical tensor language over} $\mathit{Lit}_+$.
Free and bound indices, boundaries, $\alpha$-equivalence, reparameterizations, type symbols, valencies and balanced type symbols are defined exactly as in the  intuitionistic case.
\ed
The set of intuitionistic tensor formulas is identified with a subset of classical ones by means of the familiar translation $A\multimap B=\overline{A}\parr B$.
\bd
A {\it classical tensor context} is a finite set $\Gamma$ of  tensor formulas satisfying the same requirement (\ref{indices in context}) as in the intuitionistic case. As in the intuitionistic case
  we denote $FI^\bullet(\Gamma)=\bigcup\limits_{A\in\Gamma} FI^\bullet(A)$,
$FI_\bullet(\Gamma)=\bigcup\limits_{A\in\Gamma} FI_\bullet(A)$.
A   {\it classical tensor sequent}   $\Sigma$ is an expression of the form $\vdash \Gamma$, where $\Gamma$ is a tensor context satisfying $FI^\bullet(\Gamma)=FI_\bullet(\Gamma)$.
    \ed
    Just as in the intuitionistic case, the above definition means that every index occurring in a classical tensor sequent has exactly one upper and one lower occurrence.
 Intuitionistic tensor sequents are identified with classical ones by translating $\Gamma\vdash B$ to $\vdash \overline{\Gamma}, B$, where the negated context  $\overline\Gamma$ is defined by
 $\overline{A_{(1)},\ldots,A_{(n)}}= \overline{A_{(1)}},\ldots,\overline{A_{(n)}}$. Arguably, the classical format allows more symmetric notation and compact definitions. The sequent calculus is also more compact.
\bd
The system of {\it classical  tensor term logic} ({\bf TTL}) is given by the rules in Figure \ref{TTL_rules}, where  it is assumed  that all expressions are well-formed.
  \ed
  Observe that there is the special instance $\vdash \bot^j_i,{\bf 1}^i_j$ of the $({\rm{Id}})$ axiom, which is the rule introducing the unit ${\bf 1}$.
\begin{figure}
  \centering
  \subfloat[Classical {\bf TTL} Language
\label{TTL_lang}]
{
$\begin{array}{c}
 \quad\quad\quad\quad\quad\quad \mathit{Lit}_-=\{\overline{p}|~p\in \mathit{Lit}_+\}, \quad\overline{\overline{p}}=p\mbox{ for }p\in \mathit{Lit}_+, \quad \mathit{Lit}=\mathit{Lit}_+\cup \mathit{Lit}_-.\\[.1cm]
 \quad\quad\quad\quad\quad\quad\mathit{Prop}=\mathit{Lit}\cup\{{1},\bot\},\quad v({1})=v(\bot)=(1,1),\quad \overline{{1}}=\bot,~\overline{\bot}={1}.\\[.1cm]
%\quad\quad\quad\quad\quad\quad\widetilde{\mathit{At}}=\{p_{j_1\ldots j_n}^{i_1\ldots i_m}|~p\in \mathit{Prop},v(p)=(m,n),i_1,\ldots, i_m,j_1,\ldots,j_n\in\mathit{Ind}\}.\\[.1cm]
\quad\quad\quad\quad\quad
\widetilde{\mathit{At}}=\{p_{j_1\ldots j_n}^{i_1\ldots i_m}|~p\in \mathit{Prop},i_1,\ldots, i_m,j_1,\ldots, j_n\subset\mathit{Ind},v(p)=(m,n)\}.\\[.1cm]
\widetilde{\mathit{Fm}}::=\widetilde{\mathit{At}}|(\widetilde{\mathit{Fm}}\otimes\widetilde{\mathit{Fm}})|
(\widetilde{\mathit{Fm}}\parr\widetilde{\mathit{Fm}}).\\[.1cm]
\quad\quad\quad\quad\quad\quad\overline{p^{i_1,\ldots,i_n}_{j_1,\ldots,j_m}}=\overline{p}_{{i_n},\ldots, {i_1}}^{{j_m},\ldots, {j_1}}\mbox{ for }p\in\mathit{Prop},~
%\\[.1cm]
 %\quad\quad\quad\quad\quad\quad
 \overline{A\otimes B}=\overline{B}\wp\overline{A},~
\overline{A\wp B}=\overline{B}\otimes\overline{A}.
%\overline{{1}_i^j}= \bot^i_j\,\quad
%\overline{\bot_i^jA}= {1}^i_j.
\end{array}
$
}\\
%%%%%%%%%%%%%%%%%%%%%%%%%%%%%%%%%%%%%
\subfloat[Sequent rules
  \label{TTL_rules}
]
{
$
\begin{array}{c}
\begin{array}{rl}
\cfrac{A\in At}{\vdash\overline A, A}{\rm{(Id)}}
&
\cfrac{\vdash\Gamma, A\quad \vdash\overline A,\Delta}{\vdash\Gamma,\Delta}({\rm{Cut}})
\\
\cfrac{\vdash\Gamma, A\quad \vdash B,\Delta}{\vdash\Gamma,A\otimes B,\Delta}(\otimes)
&
\cfrac{\vdash\Gamma,A, B}{\vdash\Gamma, A\parr B}(\parr)
\\
\cfrac{\vdash\Gamma, A~ j\in FI^\bullet(A)}{\vdash\Gamma,\bot^j_i,A^{[\nicefrac{i}{j}]}}({\bot}^\to)
&
\cfrac{\vdash\Gamma,A,\bot_i^j~ j\in FI_\bullet(A)}{\vdash\Gamma,A_{[\nicefrac{i}{j}]}}({\bot}_\leftarrow)
\end{array}\\[.6cm]
\mbox{In }({\bot}^\to),({\bot}_\leftarrow)~
i\not\in FI^\bullet(A)\cup FI_\bullet(A)
\end{array}
$
}
\\
%%%%%%%%%%%%%%%%%%%%%%%%%%%%%%%%%
\subfloat[Admissible rules
  \label{TTL_admissible rules}
]
{
$
\begin{array}{c}
\begin{array}{rl}
\cfrac{\vdash\Gamma, A~ j\in FI_\bullet(A)}{\vdash\Gamma,{\bot}^i_j,A_{[\nicefrac{i}{j}]}}({\bot}_\to)
&
\cfrac{\vdash\Gamma,A,{\bot}^i_j~ j\in FI^\bullet(A)}{\vdash\Gamma,A^{[\nicefrac{i}{j}]}}({\bot}^\leftarrow)
\end{array}\\[.6cm]
\mbox{In both rules above }
i\not\in FI^\bullet(A)\cup FI_\bullet(A)
\end{array}
$
}
\caption{Classical {\bf TTL}}
\end{figure}

In general, the classical system is much more convenient for
 proof-theoretic analysis (cut-elimination, proof-search etc),
 and we will use it whenever such an analysis is important. The embedding of the intuitionistic system is conservative (as we will show shortly), so the results obtained for the classical version apply to the intuitionistic fragment as well. We will sometimes use the two-sided format for {\bf TTL} as well and write $\Gamma\vdash \Delta$ for $\vdash \overline{\Gamma}, \Delta$.

\subsection{Explanation}
The explanation of the classical system is as follows. In the ordinary classical linear logic the negation or duality of formulas can be defined by
$\overline{A}\cong A\multimap \bbot$, where $\bbot$, the {\it dualizing formula} and {\it comonoidal counit}, is the dual of the monoidal unit $\mathbb{1}$,  $\bbot=\overline{\mathbb{1}}$, and it is postulated that $(A\multimap\bbot)\multimap\bbot\cong A$ for any formula $A$. (Our constants ${\bf 1}$, $\bot$ and the corresponding formulas ${\bf 1}^i_j$, $\bot^j_i$ should not be confused with the above $\mathbb{1},\bbot$, see Remark \ref{unit_remark}).

Accordingly, for interpretation of classical tensor formulas we need to choose the {\it dualizng type} $\bbot$, which  is  some set of tensor terms {\it with no indices} (this means that every term in $\bbot$ is either 1 or purely singular, see Definition \ref{tensor terms}). Classical tensor type $A$ is then a  set $A$ of terms which  is equal to its {double dual}, $A=\overline{\overline A}$, where $\overline A=A\multimap\bbot$ and $A\multimap X=\{t|~\forall s\in A~ ts\in X\}$.

The operations on classical types are defined by
$$A\parr B=\overline A\multimap B,\quad A\otimes B=\overline{\overline A\parr\overline B},\quad {\bf 1}^i_j=\overline{\overline{\{\delta^i_j\}}},\quad\bot^i_j=\overline{{\bf 1}^j_i}.$$
The counit type $\bot^i_j$ can be described explicitly as
$$\bot^i_j=\bigg\{[w_1]\cdots[w_n][w]^i_j\bigg|~[w_1]\cdots[w_n][w]\in\bbot\bigg\}.$$

 We interpret classical tensor formulas as classical tensor types with connectives corresponding to operations on types in the obvious way and with the convention that if $||A||$ is the interpretation of a tensor formula $A$ then for any $t\in||A||$ we have $FI^\bullet(t)=FI^\bullet(A)$, $FI_\bullet(t)=FI_\bullet(A)$ and that for $i\in FI^\bullet(A)$, respectively $i\in FI_\bullet(A)$, and $j$ not occurring in $A$ it holds that $t\in||A||$ iff $t^{[\nicefrac{j}{i}]}\in||A^{[\nicefrac{j}{i}]}||$, respectively
 $t_{[\nicefrac{i'}{i}]}\in||A_{[\nicefrac{i'}{i}]}||$.
 Note that it must be that indices in dual formulas $A$, $\overline A$ should match each other, $FI^\bullet(A)=FI_\bullet(\overline A)$, $FI_\bullet(A)=FI^\bullet(\overline A)$, because for $t\in ||A||$, $s\in||\overline A||$ the term $ts\in\bbot$ must have no free indices. This explains our convention for valencies of dual propositional symbols and duality of atomic formulas.

This semantics (which is actually complete) will be discussed in detail in subsequent sections.

\subsection{Simplest properties}
  As in the intuitionistic case we define {\it similarity transformation} of sequents as simultaneous replacement of all  free occurrences of an index with a fresh one, and {\it similar} sequents as those related by similarity transformations. {\it $\alpha$-Equivalence of sequents} is also defined the same as in the intuitionistic case. Then
  we list simplest properties of {\bf TTL}, parallel to those of {\bf ITTL}. Proofs are analogous to Proposition \ref{ITTL properties}.
  \bp\label{TTL properties}
  In {\bf TTL}
  \begin{enumerate}[(i)]
%  \item In a derivable sequent every free index has exactly one upper and one lower free occurrence.
   \item Similar sequents are cut-free derivable from each other;
 \item rules in Figure \ref{TTL_admissible rules} are admissible and can be emulated without use of the Cut rule;
    %\item for any formula $F$ the sequent $\vdash\overline{F},F$ is derivable;
      \item $\alpha$-equivalent formulas are provably equivalent, i.e. if $F\equiv_\alpha F'$ then the sequents $\vdash\overline F, F'$, $\vdash\overline{F'}, F$ are derivable;
          \item $\alpha$-equivalent sequents are derivable from each other;
  \item
 if   $j\in FI^\bullet(A)$, and $i$ is an index not occurring in $A$ then
  there are provable equivalences
  ${\bf 1}^i_j\otimes A \cong {\bot}^i_j\parr A\cong A^{[\nicefrac{i}{j}]}$.

  Similarly if   $j\in FI_\bullet(A)$, and $i$ is an index not occurring in $A$ then
  there are provable equivalences
  $A\otimes {\bf 1}^j_i cong  A\parr {\bot}^j_i\cong A_{[\nicefrac{i}{j}]}$. $\Box$
  \end{enumerate}
  \ep

\section{Cut-elimination}\label{cut-elimination section}
In this section we prove cut-elimination for the classical {\bf TTL}, which also implies cut-elimination  for the intuitionistic version. It will be convenient to use a number of notational and terminological conventions.

If $\Sigma$ is a tensor sequent of the form $\vdash\Gamma$ we denote $FI(\Sigma)=FI^\bullet(\Gamma)=FI_\bullet(\Gamma)$. For a context $\Gamma$ we denote $|FI(\Gamma)|=FI^\bullet(\Gamma)\cup FI_\bullet(\Gamma)$.
%If $i_1,\ldots,i_k\in FI(\Sigma)$ and $\Sigma'$ is obtained from $\Sigma$ by a chain of similarity transformations replacing $i_1,\ldots,i_k$ with some %$i_1',\ldots,i_k'\in\mathit{Ind}$ respectively, we denote $\Sigma'=\Sigma[\nicefrac{i_1'}{i_1},\ldots,\nicefrac{i_k'}{i_k}]$.
%Finally, if $FI(\Sigma)=\{i_1,\ldots,i_n\}$, the map $\phi:FI(\Sigma)\to FI(\Sigma')$ is a bijection and $\Sigma'=\Sigma[\nicefrac{\phi(i_1)}{i_1},\ldots,\nicefrac{\phi(i_n)}{i_n}]$  we write simply $\Sigma'=\Sigma[\phi]$.

We treat sequent calculus derivations as rooted trees whose nodes are labeled with sequents, with the root corresponding to the conclusion and leaves to axioms.
 Given a derivation $\pi$, we say that the index $i$ occurs, respectively occurs freely, in $\pi$ if it occurs, respectively occurs freely, in the sequent labeling some node of $\pi$. We denote the set of indices occurring in $\pi$ freely as $FI(\pi)$. If $i,i'\in \mathit{Ind}$ we denote the {labeled tree obtained from $\pi$ by replacing all free occurrences of $i$ with $i'$} as $\pi[\nicefrac{i'}/{i}]$. Similarly, we denote the expression
obtained from a sequent $\Sigma$ by replacing all free occurrences of $i$ with $i'$ as $\Sigma[\nicefrac{i'}{i}]$.

We use the following notation for iterated substitutions. Capital Latin letters ($I,J,\ldots$) stand for sequences of pairwise distinct indices and the expression $[\nicefrac{I'}{I}]$, where $I=\{i_1,\ldots,i_n\}$, $I=\{i_1',\ldots,i_n'\}$, stands for ${[\nicefrac{i_1'}{i_1}],\ldots[\nicefrac{i_n'}{i_n}]}$. We use the same notation when $I',I$ are (unordered) index sets, rather than sequences, but the bijection $i_\alpha\mapsto i_\alpha'$, $\alpha=1,\ldots,n$, is clear from the context.

\subsection{Alternative formulation}
\begin{figure}
\centering
\subfloat
[
${\bf TTL}'$ rules
\label{primed system}
]
{
$
\begin{array}{c}
\cfrac
    {
    \vdash \Gamma,A, B
    }
    {
    \vdash\Gamma, A^{[\nicefrac{K}{I}]}_{[\nicefrac{L}{J}]}\parr B_{[\nicefrac{K}{I}]}^{[\nicefrac{L}{J}]}
    }
    (\parr_{\equiv\alpha})
~~
\cfrac
    {
    \vdash\Gamma,A~\vdash B,\Theta
    }
    {
    \vdash\Gamma, A^{[\nicefrac{K}{I}]}_{[\nicefrac{L}{J}]}\otimes B_{[\nicefrac{K}{I}]}^{[\nicefrac{L}{J}]},\Theta
    }
    (\otimes_{\equiv\alpha})\\[.5cm]
   \mbox{where }I\subseteq FI^\bullet(A)\cap FI_\bullet(B), J\subseteq FI_\bullet(A)\cap FI^\bullet(B),
   \\
   K\cap( FI^\bullet(A,B)\cup FI_\bullet(A,B))=L\cap( FI^\bullet(A,B)\cup FI_\bullet(A,B))=\emptyset.
   \end{array}
$
}
\\[.3cm]
%%%%%%%%%%%%%%%%%%%%%%%%%%%%%%%%%%%%%%%%%%%%
\subfloat[Emulating the $(\otimes_{\equiv\alpha})$ rule in {\bf TTL}
\label{primed system emulation figure}]
{
$\begin{array}{c}
k\in FI^\bullet(\Gamma)\cap FI_\bullet(\Gamma),
l\in FI_\bullet(\Delta)\cap FI^\bullet(\Delta), k',l'\mbox{ fresh},\\
i\in FI_\bullet(\Gamma)\cap FI^\bullet(A)\cap FI_\bullet(B)\cap FI^\bullet(\Delta),
\\
j\in FI^\bullet(\Gamma)\cap FI_\bullet(A)\cap FI^\bullet(B)\cap FI_\bullet(\Delta):\\[.2cm]
\def\fCenter{\ \vdash\ }
\def\labelSpacing{.5pt}
\def\ScoreOverhang{.1pt}
\Axiom$\fCenter \Gamma,A$
\RightLabel{($\sim$)}
\UnaryInf$ \fCenter\Gamma_{[\nicefrac{k'}{k}]}^{[\nicefrac{k'}{k}]},A$
\RightLabel{($\sim$)}
\UnaryInf$ \fCenter(\Gamma_{[\nicefrac{k'}{k}]}^{[\nicefrac{k'}{k}]})_{[\nicefrac{k}{i}]}^{[\nicefrac{l}{j}]},A^{[\nicefrac{k}{i}]}_{[\nicefrac{l}{j}]}$
\Axiom$\fCenter B,\Delta$
\RightLabel{($\sim$)}
\UnaryInf$\fCenter B,\Delta^{[\nicefrac{l'}{l}]}_{[\nicefrac{l'}{l}]}$
\RightLabel{($\sim$)}
\UnaryInf$ \fCenter B_{[\nicefrac{k}{i}]}^{[\nicefrac{l}{j}]},
(\Delta^{[\nicefrac{l'}{l}]}_{[\nicefrac{l'}{l}]})^{[\nicefrac{k}{i}]}_{[\nicefrac{l}{j}]}$
\RightLabel{($\otimes$)}
\BinaryInf$\fCenter(\Gamma_{[\nicefrac{k'}{k}]}^{[\nicefrac{k'}{k}]})_{[\nicefrac{k}{i}]}^{[\nicefrac{l}{j}]},
A^{[\nicefrac{k}{i}]}_{[\nicefrac{l}{j}]}
\otimes B_{[\nicefrac{k}{i}]}^{[\nicefrac{l}{j}]},(\Delta^{[\nicefrac{l'}{l}]}_{[\nicefrac{l'}{l}]})^{[\nicefrac{k}{i}]}_{[\nicefrac{l}{j}]}
$
\RightLabel{($\sim$)}
\UnaryInf$ \fCenter\Gamma^{[\nicefrac{k'}{k}]}_{[\nicefrac{k'}{k}]},A^{[\nicefrac{k}{i}]}_{[\nicefrac{l}{j}]}
\otimes B_{[\nicefrac{k}{i}]}^{[\nicefrac{l}{j}]},\Delta^{[\nicefrac{l'}{l}]}_{[\nicefrac{l'}{l}]}$
\RightLabel{($\sim$)}
\UnaryInf$ \fCenter\Gamma,A^{[\nicefrac{k}{i}]}_{[\nicefrac{l}{j}]}
\otimes B_{[\nicefrac{k}{i}]}^{[\nicefrac{l}{j}]},\Delta$
\DisplayProof
\end{array}
$
}
\caption{Alternative formulation of {\bf TTL}}
\end{figure}
For proving cut-elimination it is convenient to give an alternative formulation for (classical) {\bf TTL},
where the $(\otimes)$ and the $(\parr)$ rules are replaced with their closures under $\alpha$-equivalence of sequents.
\bd
The system ${\bf TTL}'$  is obtained from  {\bf TTL}  by  replacing the $(\otimes)$ and the $(\parr)$ rules with their  closures $(\otimes_{\equiv_\alpha})$, $(\parr_{\equiv_\alpha})$ under $\alpha$-equivalence shown in Figure \ref{primed system}.
\ed
\bp\label{TTL'->TTL}
A sequent is derivable, respectively cut-free derivable, in ${\bf TTL}'$ iff it is derivable, respectively cut-free derivable, in ${\bf TTL}$.
\ep
{\bf Proof} {\bf TTL} rules  are special instances  of  ${\bf TTL}'$ rules, so the direction from ${\bf TTL}$ to ${\bf TTL}'$ is trivial. For the other direction it can be observed that the $(\otimes)'$, $(\parr)'$ rules can be cut-free emulated in {\bf TTL} using Proposition \ref{TTL properties}, clause {(i)}.  In Figure \ref{primed system emulation figure} we show
a generic example for the $(\otimes_{\equiv\alpha})$ rule, where, for notational simplicity, we assume that all sequences $K,L,I,J$ are singletons $k,l,i,j$. Note, however, that we should take into account that the indices $k$ and $l$ may occur in $\Gamma$  and $\Delta$ respectively.
 The indices $k'$ and $l'$ in Figure \ref{primed system emulation figure} are fresh; they are needed in order to avoid possible index collisions. $\Box$

The advantage of using ${\bf TTL}'$ is that derivations in this system are better behaving under index substitutions.
\bp\label{substituion in derivation}
If $\pi$ is a ${\bf TTL}'$ derivation, $i\in FI(\pi)$ and $i'\in\mathit{Ind}$ does not occur in $\pi$ then $\pi[\nicefrac{i'}{i}]$ is a valid ${\bf TTL}'$ derivation.
\ep
{\bf Proof} Induction on $\pi$.
If $\pi$ is an axiom the statement is obvious.

Assume that $\pi$
was obtained from derivations $\pi_1,\pi_2$ with conclusions $\Sigma,\Pi$ of the form $\vdash\Gamma,A$, $\vdash B,\Delta$ respectively by means of the $(\otimes_{\equiv\alpha})$ rule.
Then the conclusion $\Lambda$ of $\pi$ has the form
$\vdash\Gamma,F,\Delta$, where $F=A^{[\nicefrac{K}{I}]}_{[\nicefrac{L}{J}]}\otimes B_{[\nicefrac{K}{I}]}^{[\nicefrac{L}{J}]}$ for some
$I\subseteq FI^\bullet(A)\cap FI_\bullet(B)$, $J\subseteq FI_\bullet(A)\cap FI^\bullet(B)$, and $K,L$ such that   $K\cap |FI(A,B)|=L\cap |FI(A,B)|=\emptyset$,

By the induction hypothesis we have valid ${\bf TTL}'$ derivations $\pi_1'=\pi_1[\nicefrac{i'}{i}]$, $\pi_2'=\pi_2[\nicefrac{i'}{i}]$ with conclusions $\Sigma'=\Sigma[\nicefrac{i'}{i}]$, $\Pi'=\Pi[\nicefrac{i'}{i}]$ respectively. The statement will be proven if we show that the sequent $\Lambda'=\Lambda[\nicefrac{i'}{i}]$ is derivable from $\Sigma'$ and $\Pi'$ by the $(\otimes_{\equiv\alpha})$ rule.

The most involved case is when $i$ has free occurrences both in $\Sigma$ and $\Pi$. In this case there must be exactly one free occurrence of $i$ in each of the four contexts $\Gamma$, $A$, $B$ and $\Delta$, otherwise $\Lambda$ is not well-formed. Assume for definiteness that $i\in FI_\bullet(\Gamma)$, hence
$i\in FI^\bullet(A)$, $i\in FI_\bullet(B)$, $i\in FI^\bullet(\Delta)$. Then $\Sigma',\Pi',\Lambda'$ are respectively of the forms $\vdash\Gamma_{[\nicefrac{i'}{i}]},A'$, $\vdash B',\Delta^{[\nicefrac{i'}{i}]}$ and $\vdash\Gamma_{[\nicefrac{i'}{i}]},F,\Delta^{[\nicefrac{i'}{i}]}$,
where $A'=A^{[\nicefrac{i'}{i}]}$, $B'=B_{[\nicefrac{i'}{i}]}$. In particular $i'\in FI^\bullet(A')\cap FI_\bullet(B')$ and
$F=A'^{[\nicefrac{i}{i'}]}\otimes B'_{[\nicefrac{i}{i'}]}$. It follows that
 indeed $\Lambda'$ is derivable from $\Sigma'$ and $\Pi'$ by the $(\otimes_{\equiv\alpha})$ rule. Other cases are straightforward.

Other rules are treated similarly.
$\Box$

\subsection{Cut-elimination in ${\bf TTL}'$}
We define the {\it size of a derivation} in ${\bf TTL}'$ as follows. If the derivation $\pi$ is an axiom then the size $\mathit{size}(\pi)$ of $\pi$ is $1$. If $\pi$ is obtained from a derivation $\pi'$ by a single-premise rule then $\mathit{size}(\pi)=\mathit{size}(\pi')+1$. If $\pi$ is obtained from derivations $\pi_1$, $\pi_2$ by a two-premise rule then
$\mathit{size}(\pi)=\mathit{size}(\pi_1)+\mathit{size}(\pi_2)+1$.

%Note that the $(\equiv_\beta)$ rules do not count. (Essentially, in the analysis below we identify $\beta$-equivalent tensor terms or sequents.)

Given  a ${\bf TTL}'$ derivation with two subderivations $\pi_1$, $\pi_2$ of the sequents
$\vdash \Gamma,A$ and $\vdash \overline A, \Delta$ respectively followed by  the Cut rule, we say that this application of the Cut rule   is a {\it first type side cut} if the last rule in $\pi_1$ or $\pi_2$  does not involve the cut-formula. The application is a {\it second type side cut} if
the last rule in $\pi_1$ or $\pi_2$ is the $(\bot^\to)$ or the $(\bot_\leftarrow)$ rule modifying the cut-formula.
  An application of the Cut rule that is not a side cut is a {\it principal  cut}.
\bl\label{side cuts}
There is an algorithm transforming a ${\bf TTL}'$ derivation $\pi$ of a sequent $\vdash\Gamma$ to a derivation $\pi'$ of  the same sequent such that $\mathit{size}(\pi')=\mathit{size}(\pi)$ and all cuts in $\pi'$ principal.
\el
Let $\pi_1,\ldots,\pi_n$ be all subderivations of $\pi$ ending with side cuts and $N=\sum\mathit{size}(\pi_i)$. We will use $N$ as the induction parameter for proving termination of the algorithm.

The algorithm consists in permuting a side cut with the preceding rule, which may be accompanied with tours of index renaming.

Cut-elimination step for the first type side cut is shown in Figure \ref{side cut figure}, where we have an application of the $(\parr_{\equiv\alpha})$ rule to the sequent $\vdash\Gamma,F,G,A$ followed by a side cut with the sequent $\vdash \overline A,\Delta$. In general, we cannot cut the premise $\vdash\Gamma, F,G,A$ of the $(\parr_{\equiv_\alpha})$ rule with the right premise $\vdash\overline A,\Delta$ of the side cut, because this might  introduce forbidden repetitions of indices (if $|FI(\Gamma,\Delta)|\cap(I\cup J)\not=\emptyset$). Therefore we use Proposition \ref{substituion in derivation} and at first replace all free occurrences of $I,J$ in the left subderivation with fresh indices.
The induction parameter $N$ for the new derivation is smaller at least by 1. Other instances of the first type cut are treated similarly or easier.

Cut-elimination step for  the second  type is shown in Figure \ref{side cut figure2}, where  we have an application of the $(\bot_\leftarrow)$ rule to the sequent $\vdash\overline{A},\bot^j_i,\Delta$
 on the  right, which is followed by a side cut with the sequent $\vdash  \Gamma,A^{[\nicefrac{i}{j}]}$.
Again, after the transformation the induction parameter $N$ is smaller at least by 1. The case of the $(\bot^\to)$ rule is treated similarly.
$\Box$
      \begin{figure}
      \centering
       \subfloat[{Transforming a 1st type side cut}
\label{side cut figure}
]
{
      $
   \begin{array}{c}
   I\subseteq FI^\bullet(F)\cap FI_\bullet(G),~J\subseteq FI_\bullet(F)\cap FI^\bullet(G),\\
   K\cap |FI(A,B)|=L\cap |FI(A,B)|=\emptyset,~
   K',L'\mbox{ fresh}:
  \\[.3cm]
%(\nu,\mu)\in FI(\sigma),~ (\mu,\nu)\in FI(B),~\\
%\def\ScoreOverhang{.1pt}
%\def\fCenter{\ \ \ }
\AxiomC{$\pi_1$}
%\noLine
%%%%\def\fCenter{\ \vdash\ }
\UnaryInfC{$\vdash\Gamma,F,G,A$}
\RightLabel{($\parr_{\equiv\alpha}$)}
\UnaryInfC{$ \vdash\Gamma,F^{[\nicefrac{K}{I}]}_{[\nicefrac{L}{J}]}\parr G_{[\nicefrac{K}{I}]}^{[\nicefrac{L}{J}]},A$}
\AxiomC{$\pi_2$}
%\noLine
%%\def\fCenter{\ \vdash\ }
\UnaryInfC{$\vdash\overline A,\Delta$}
\RightLabel{(${\rm{Cut}}$)}
\insertBetweenHyps{\hskip 1pt}
\BinaryInfC{$ \vdash \Gamma,F^{[\nicefrac{K}{I}]}_{[\nicefrac{L}{J}]}\parr G_{[\nicefrac{K}{I}]}^{[\nicefrac{L}{J}]},\Delta$}
\DisplayProof
\Longrightarrow
\def\ScoreOverhang{.1pt}
\AxiomC{$\pi_1[\nicefrac{K'}{K},\nicefrac{L'}{L}]$}
%\noLine
%\def\fCenter{\ \vdash\ }
\UnaryInfC{$\vdash\Gamma,F^{[\nicefrac{K'}{I}]}_{[\nicefrac{L'}{J}]},G_{[\nicefrac{K'}{I}]}^{[\nicefrac{L'}{J}]},A$}
\AxiomC{$\pi_2$}
%\noLine
%\def\fCenter{\ \vdash\ }
\UnaryInfC{$\vdash \overline A,\Delta$}
\RightLabel{(${\rm{Cut}}$)}
\insertBetweenHyps{\hskip 5pt}
\BinaryInfC{$ \vdash \Gamma,F^{[\nicefrac{K'}{I}]}_{[\nicefrac{L'}{J}]},G_{[\nicefrac{K'}{I}]}^{[\nicefrac{L'}{J}]},\Delta$}
\RightLabel{($\parr_{\equiv\alpha}$)}
\UnaryInfC{$  \vdash\Gamma,F^{[\nicefrac{K}{I}]}_{[\nicefrac{L}{J}]}\fCenter\parr G_{[\nicefrac{K}{I}]}^{[\nicefrac{L}{J}]},\Delta$}
\DisplayProof
\end{array}
   $
   }
   \\[.3cm]
%%%%%%%%%%%%%%%%%%%%%%%%%%%%%%%%%%%%%%%%
\subfloat[{Transforming a 2nd type side cut}
\label{side cut figure2}
]
{
   $
   \begin{array}{c}
   j\in FI^\bot(A),~i\not\in|FI(A)|,~j'\mbox{ fresh}:
    \\[.3cm]
 \def\ScoreOverhang{.1pt}
\Axiom$\fCenter\pi_1$
%\noLine
\def\fCenter{\ \vdash\ }
\UnaryInf$\fCenter \Gamma,A^{[\nicefrac{i}{j}]}$
\def\fCenter{\ \ \ }
\Axiom$\fCenter\pi_2$
%\noLine
\def\fCenter{\ \vdash\ }
\UnaryInf$\fCenter\overline{A},\bot^j_i,\Delta$
\RightLabel{($\bot_\leftarrow$)}
\UnaryInf$ \fCenter {\overline{A}_{[\nicefrac{i}{j}]}},\Delta$
\RightLabel{(${\rm{Cut}}$)}
\insertBetweenHyps{\hskip 5pt}
\BinaryInf$ \fCenter \Gamma,\Delta$
\DisplayProof
\Longrightarrow
\def\ScoreOverhang{.1pt}
\def\fCenter{\ \ \ }
\Axiom$\fCenter\pi_1[\nicefrac{j'}{i}]$
%\noLine
\def\fCenter{\ \vdash\ }
\UnaryInf$\fCenter \Gamma_{[\nicefrac{j'}{i}]},A^{[\nicefrac{j'}{j}]}$
\def\fCenter{\ \ \ }
\Axiom$\fCenter\pi_2{[\nicefrac{j'}{j}]}$
%\noLine
\def\fCenter{\ \vdash\ }
\UnaryInf$\fCenter \overline A_{[\nicefrac{j'}{j}]},\bot^{j'}_i,\Delta$
\RightLabel{(${\rm{Cut}}$)}
\insertBetweenHyps{\hskip 5pt}
\BinaryInf$ \fCenter \Gamma_{[\nicefrac{j'}{i}]},\bot^{j'}_i,\Delta$
\RightLabel{($\bot_\leftarrow$)}
\UnaryInf$ \fCenter \Gamma,\Delta$
\DisplayProof
\end{array}
   $
   }   %
  \\[.3cm]
%%%%%%%%%%%%%%%%%%%%%%%%%%%%%%%%%%%%%%%%
  \subfloat[{Transforming a principal cut}
\label{principal cut figure}
]
{
   $
   \begin{array}{c}
   I\subseteq FI^\bullet(A)\cap FI_\bullet(B),~J\subseteq FI_\bullet(A)\cap FI^\bullet(B),\\
   K\cap |FI(A,B)|=L\cap |FI(A,B)|=K'\cap |FI(A,B)|=L'\cap |FI(A,B)|=\emptyset:
  \\[.3cm]
%(\nu,\mu)\in FI(\sigma),~ (\mu,\nu)\in FI(B),~\\
\def\ScoreOverhang{.1pt}
\def\fCenter{\ \ \ }
\Axiom$\fCenter\pi_1$
\def\fCenter{\ \vdash\ }
\UnaryInf$\fCenter\Gamma,A^{[\nicefrac{K}{I}]}_{[\nicefrac{L}{J}]}$
\def\fCenter{\ \ \ }
\Axiom$\fCenter\pi_2$
\def\fCenter{\ \vdash\ }
\UnaryInf$\fCenter \Delta,B_{[\nicefrac{K}{I}]}^{[\nicefrac{L}{J}]}$
\RightLabel{$(\otimes_{\equiv\alpha})$}
\BinaryInf$ \fCenter\Gamma,\Delta,A\otimes B$
\def\fCenter{\ \ \ }
\Axiom$\fCenter\pi_3$
\def\fCenter{\ \vdash\ }
\UnaryInf$\fCenter \overline{A}_{[\nicefrac{K'}{I}]}^{[\nicefrac{L'}{J'}]},\overline{B}^{[\nicefrac{K'}{I}]}_{[\nicefrac{L'}{J'}]},\Theta$
\RightLabel{$(\parr_{\equiv\alpha})$}
\UnaryInf$\fCenter\overline A\parr\overline B,\Theta$
\RightLabel{(${\rm{Cut}}$)}
\insertBetweenHyps{\hskip 1pt}
\BinaryInf$ \fCenter \Gamma,\Delta,\Theta$
\DisplayProof
\\
\Longrightarrow
\\
\def\ScoreOverhang{.1pt}
\def\fCenter{\ \ \ }
\Axiom$\fCenter\pi_1[\nicefrac{K'}{K},\nicefrac{L'}{L}]$
\def\fCenter{\ \vdash\ }
\UnaryInf$\fCenter\Gamma_{[\nicefrac{K'}{K}]}^{[\nicefrac{L'}{L}]},A^{[\nicefrac{K'}{I}]}_{[\nicefrac{L'}{J}]}$
\def\fCenter{\ \ \ }
\Axiom$\fCenter\pi_3$
\def\fCenter{\ \vdash\ }
\UnaryInf$\fCenter \overline{A}_{[\nicefrac{K'}{I}]}^{[\nicefrac{L'}{J}]},\overline{B}^{[\nicefrac{K'}{I}]}_{[\nicefrac{L'}{J'}]},\Theta$
\RightLabel{$({\rm{Cut}})$}
\insertBetweenHyps{\hskip 1pt}
\BinaryInf$\fCenter\Gamma_{[\nicefrac{K'}{K}]}^{[\nicefrac{L'}{L}]},\overline{B}^{[\nicefrac{K'}{I}]}_{[\nicefrac{L'}{J}]},\Theta$
\def\fCenter{\ \ \ }
\Axiom$\fCenter\pi_2[\nicefrac{K'}{K},\nicefrac{L'}{L}]$
\def\fCenter{\ \vdash\ }
\UnaryInf$\fCenter \Delta^{[\nicefrac{K'}{K}]}_{[\nicefrac{L'}{L}]},B_{[\nicefrac{K'}{I}]}^{[\nicefrac{L'}{J}]}$
\RightLabel{(${\rm{Cut}}$)}
\insertBetweenHyps{\hskip 1pt}
\BinaryInf$ \fCenter \Gamma_{[\nicefrac{K'}{K}]}^{[\nicefrac{L'}{L}]},\Delta^{[\nicefrac{K'}{K}]}_{[\nicefrac{L'}{L}]},\Theta$
\DisplayProof
\end{array}
   $
   } \caption{Cut-elimination steps}
%\label{side cut figure}
 \end{figure}
 \bl\label{principal cuts}
  There is an algorithm transforming a  ${\bf TTL}'$ derivation $\pi$  with cuts into a  ${\bf TTL}'$ derivation $\pi'$ such that the conclusions of $\pi$ and $\pi'$ are similar and $\mathit{size}(\pi')<\mathit{size}(\pi)$, provided  that all cuts in $\pi$ are principal.
 \el
{\bf Proof} The algorithm consists in pushing principal cuts upwards to axioms after appropriate index renaming.

If the application of the cut rule is principal and the cut-formulas are atomic, then at least   one of the premises is an axiom and it can be erased.
If the cut-formulas are compound it must be that they have been introduced right above the cut then the rules above the cut  by the $(\otimes_{\equiv_\alpha})$ and
$(\parr_{\equiv\alpha})$ rules.
This case is shown in Figure \ref{principal cut figure}. $\Box$
\bc
For any ${\bf TTL}'$ derivation  $\pi$ there exists a cut-free derivation whose conclusion is similar to that of $\pi$. $\Box$
\ec
\bc
The system ${\bf TTL}'$ is cut-free.
\ec
{\bf Proof} Observe that Proposition \ref{TTL properties}, clause (i) (that similar sequents are cut-free derivable from each other) holds for ${\bf TTL}'$. $\Box$
\bc The system {\bf TTL} is cut-free.
\ec
{\bf Proof} Follows from the preceding corollary and Proposition \ref{TTL'->TTL}. $\Box$
\bc
The system {\bf ITTL} is a conservative fragment of {\bf TTL}. In particular, {\bf ITTL} is cut-free. $\Box$
\ec
In the sequel we will be also interested in derivations from {\it non-logical axioms}. For this case we have the following version of cut-elimination, which is proven identically.
\bl\label{cut-elimination with axioms}
Let $\Xi$ be a set of tensor sequents. If a sequent $\vdash\Gamma$ is {\bf TTL} derivable from $\Xi$ then there exists a {\bf TTL} derivation $\pi$ of $\vdash\Gamma$ from $\Xi$ such that for any application of the Cut rule in $\pi$ at least  one of the two premises is an element of $\Xi$ and if the cut-formula is atomic non-unit (i.e. not of the form ${\bf 1}^i_j$ or $\bot^i_j$) then both premises are elements of $\Xi$. $\Box$
\el

\section{Decidability}
Although {\bf TTL} is cut-free, it does not  satisfy the subformula property in the direct sense, unless we agree
that the counit $\bot^j_i$ is a subformula of everything, which seems rather exotic. The reason is in the reduction rule $(\bot_\leftarrow)$. This rule also makes decidability of the system rather non-obvious. Our current goal is to show that applications of this rule can be heavily restricted, which makes proof-search finite and the system decidable.
\bd
An application
$$
\cfrac{\vdash\Gamma,A,\bot^j_i\quad j\in FI_\bullet(A)}{\vdash\Gamma, A_{[\nicefrac{i}{j}]}}(\bot_\leftarrow)
$$
of the $(\bot_\leftarrow)$ rule in a {\bf TTL} proof is {\it essential} if there is a bound occurrence of the index $j$ in $\Gamma$. Otherwise it is {\it inessential}.
\ed
Our first goal is to show that inessential applications of the rule can be eliminated from the derivation.
\bp
Assume that $\vdash\Gamma$ is a  sequent
 derivable in ${\bf TTL}$ with a cut-free derivation $\pi$ and the index $j$ occurs in $\Gamma$ freely and has no bound occurrences. Let  $i$ be a fresh index.
 The sequent   $\vdash\Gamma^{[\nicefrac{i}{j}]}_{[\nicefrac{i}{j}]}$ is derivable in {\bf TTL} with a derivation $\pi'$, which is obtained from $\pi$ by replacing the index  $j$ with $i$.
\ep
{\bf Proof} Induction on $\pi$. The essential step is the $(\otimes)$ rule.

Assume that $\pi$ is obtained from derivations $\pi_1$, $\pi_2$ with conclusions $\vdash\Gamma_1$, $\vdash\Gamma_2$ respectively by the $(\otimes)$ rule. Every occurrence of $j$ in $\vdash\Gamma$ comes from one of the premises $\vdash\Gamma_1$, $\vdash\Gamma_2$. Assume for definiteness that the upper free occurrence of $j$ in $\vdash\Gamma$ comes from $\vdash\Gamma_1$. There is also a lower free occurrence of $j$ in $\vdash\Gamma_1$, and since there are no bound occurrences of $j$ in $\vdash\Gamma$ the latter occurrence remains free in $\vdash\Gamma$. It follows that both occurrences of $j$ in $\vdash\Gamma$ come from $\vdash\Gamma_1$. Also, there is no bound occurrence of $j$ in $\vdash\Gamma_1$, because otherwise it would be inherited by $\Gamma$. Thus the induction hypothesis can be applied to  $\Gamma_1$ and $\pi_1$, and we have a derivation $\pi_1'$ of $\vdash(\Gamma_1)^{[\nicefrac{i}{j}]}_{[\nicefrac{i}{j}]}$ obtained from $\pi_1$ by replacing $j$ with $i$. The desired derivation $\pi'$ of $\vdash\Gamma^{[\nicefrac{i}{j}]}_{[\nicefrac{i}{j}]}$ is obtained from $\pi_1'$ and $\pi_2$ by the $(\otimes)$ rule. $\Box$
\bp
Let a sequent $$\vdash\Gamma,A,\bot^{k_0}_{k_1},\bot^{k_1}_{k_2},\ldots,\bot^{k_{n-1}}_{k_n},$$ where $k_0\in FI_\bullet(A)$, $k_n\not\in FI^\bullet(A)\cup FI_\bullet(A)$ and there is no bound occurrence of $k_0$ in $\Gamma$, be derivable in {\bf TTL}  with a cut-free derivation $\pi$.
Then the sequent $\vdash\Gamma, A_{[\nicefrac{k_n}{k_0}]}$
is derivable in ${\bf TTL}$ with a cut-free derivation of the size smaller than  $\mathit{size}(\pi)+n-1$.
\ep
{\bf Proof} Induction on the size of $\pi$. We consider essential induction steps.

 Assume that    $\pi$ is obtained from a derivation $\pi'$ by the $(\bot^\to)$ introducing the formula $\bot^{k_0}_{k_1}$.
 This means that $\Gamma=\Gamma', B^{[\nicefrac{k_1}{k_0}]}$, and $\pi'$ has the conclusion
$$\vdash\Gamma',B,A,\bot^{k_1}_{k_2},\ldots,\bot^{k_{n-1}}_{k_n},$$ where $k_0\in FI^\bullet(B)$, $k_1\not\in FI^\bullet(B)\cup FI_\bullet(B)$.
 By the preceding proposition
we have a derivation $\pi''$ of $$\vdash\Gamma',B^{[\nicefrac{k_1}{k_0}]},A_{[\nicefrac{k_1}{k_0}]},\bot^{k_1}_{k_2},
\ldots,\bot^{k_{n-1}}_{k_n}$$ obtained from $\pi'$ by replacing $k_0$ with $k_1$. The size of $\pi''$ is strictly less than that of $\pi$. We obtain a derivation of  the desired sequent $\vdash\Gamma,A_{[\nicefrac{k_n}{k_0}]}$ by applying to $\pi''$ the $(\bot_\leftarrow)$ rule $n-1$ time.

 %If    $\pi$ is obtained from a derivation $\pi'$\ by the $(\bot^\to)$ introducing one of the formulas $\bot^{k_\alpha}_{k_{\alpha+1}}$, $\alpha>0$, the statement follows from the induction hypothesis.

 Assume that  $\pi$ is obtained from a derivation $\pi'$ by the $(\bot_\leftarrow)$ rule introducing one of the formulas $\bot^{k_\alpha}_{k_{\alpha+1}}$. For definiteness, let $\alpha=0$.
Then  the conclusion of $\pi'$ has the form $$\vdash\Gamma,A,\bot^{k_0}_{k},\bot_{k_1}^{k},\bot^{k_1}_{k_2},\ldots,\bot^{k_{n-1}}_{k_n}.$$ By the induction hypothesis there is a derivation of
$\vdash\Gamma,A_{[\nicefrac{k_n}{k_0}]}$ with size smaller than $\mathit{size}(\pi')+n=\mathit{size}(\pi)+n-1$.

Assume that $\pi$ is obtained from derivations $\pi_1,\pi_2$ by the $(\otimes)$ rule. For simplicity assume that the $(\otimes)$ rule does not modify the formula $A$ (the other case is very similar). Let $\pi_1$ be the derivation whose conclusion $\vdash\Gamma_1$ contains $A$. Since $A$ has a lower occurrence of $k_0$, it follows that $\Gamma_1$ has also an upper occurrence of $k_0$, which means that $\Gamma_1$ contains the formula $\bot^{k_0}_{k_1}$ (because there are no other occurrences of $k_0$). Thus, the conclusion $\vdash\Gamma_1$ of $\pi_1$ is of the form $\vdash\Gamma_1',A,\bot^{k_0}_{k_1}$. By the induction hypothesis we have a derivation
$\pi_1'$ with conclusion $\vdash\Gamma_1',A_{[\nicefrac{k_1}{k_0}]}$ of size smaller than that of $\pi_1$. Then, applying the $(\otimes)$ rule to $\pi_1'$ and $\pi_2$, we obtain a derivation $\pi'$ of the sequent
$\vdash\Gamma,A_{[\nicefrac{k_1}{k_0}]},\bot^{k_1}_{k_2},\ldots,\bot^{k_{n-1}}_{k_n}$ whose size is smaller than that of $\pi$. We obtain a derivation of  the desired sequent $\vdash\Gamma,A_{[\nicefrac{k_n}{k_0}]}$ by applying to $\pi''$ the $(\bot_\leftarrow)$ rule $n-1$ time.

Other cases are very easy. $\Box$
\bc
Let a sequent $\vdash\Gamma,A,\bot^j_i$, where $j\in FI_\bullet(A)$, $i\not\in FI^\bullet(A)\cup FI_\bullet(A)$ and there is no bound occurrence of $j$ in $\Gamma$, be derivable in {\bf TTL}  with a cut-free derivation $\pi$.
The sequent $\vdash\Gamma, A_{[\nicefrac{i}{j}]}$
is derivable in ${\bf TTL}$ with a cut-free derivation of the size smaller than that of $\pi$. $\Box$
\ec
\bc\label{essential reductions}
Any {\bf TTL} derivable sequent is derivable with a cut-free derivation where all applications of the $(\bot_\leftarrow)$ rule are essential. $\Box$
\ec
Thus, when searching for a {\bf TTL} derivation of a tensor sequent we may restrict to cut-free derivations where all applications of the $(\bot_\leftarrow)$ rule are essential. In such a derivation premises of any rule are built from subformulas of the conclusion and, possibly, counits $\bot^j_i$, where the indices $i,j$ occur in the conclusion. In particular there are only finitely many possible rule applications with a given conclusion. In order to guarantee that the proof-search stops we introduce a specific complexity measure of the sequent.

Let the {\it complexity} $c(A)$ of a formula $A$ be defined by
$$
c(A)=1\mbox{ for }A\mbox{ atomic},~c(A_1\otimes A_2)=c(A_1\parr A_2)=c(A_1)+c(A_2)+1.
$$
Let us say that an index $j$ is {\it suspicious} for the  sequent $\vdash\Gamma$  if there is a bound occurrence of $j$ in $\Gamma$, but the formula $\bot^j_i$ does not occur in $\Gamma$ for any $i$. Let $d(\Gamma)$ be the number of indices suspicious for $\vdash \Gamma$. Finally, let the {\it complexity} $c(\Gamma)$ of $\vdash\Gamma$ be defined as $c(\Gamma)=\sum\limits_{A\in\Gamma}c(A)+2d(\Gamma)$. Then if $\vdash\Gamma$ is a conclusion of any {\bf TTL} rule that is not a cut or an inessential application of the $(\bot_\leftarrow)$ rule, the premises of the rule have complexity strictly less than $\vdash\Gamma$. This provides us with a decision procedure for {\bf TTL} and, by conservativity, for {\bf ITTL} as well.
\bc
The systems {\bf TTL} and {\bf ITTL} are decidable. $\Box$
\ec

\section{Semantics}
We have already sketched the ``naive'' semantics for {\bf ITTL}, which is sound, but whose completeness is under question. In this section we slightly refine the above and define semantics for both {\bf TTL} and {\bf ITTL}, which is definitely complete when the language is balanced. (Completeness will be proven in the end of the paper.)

\subsection{Semantic types}
\bd
Given a terminal alphabet $T$, a {\it semantic pretype} $A$ over $T$ is a set $A$
 of tensor terms (over $T$) having common boundary, which is called the {\it boundary} of $A$. We denote the boundary of $A$ as $FI(A)$ and its two components as $FI^\bullet(A)$ and $FI_\bullet(A)$.
 \ed
 In order to define a system of semantic {\it types} that will interpret tensor formulas we need to choose a special pretype $\bbot$ with empty boundary.
\bd
 A {\it semantic type system} consists of a terminal alphabet $T$ and a pretype $\bbot$ over $T$, called the  {\it dualizing type}, satisfying  $FI^\bullet(\bbot)=FI_\bullet(\bbot)=\emptyset$.

Given a semantic type system $(T,\bbot)$ and a pretype $A$ over $T$, the {\it dual}  $\overline A$ of  $A$ is the set of terms (over $T$) defined by
$\overline A=\{s|~\forall t\in A~ts\in\bbot\}$.

The pretype $A$ is a {\it semantic type} if $A=\overline{\overline A}$.
\ed

In the following we assume that a semantic type system is given. We will omit the prefix ``semantic'', unless it leads to a confusion.
\bp\label{dual pretypes proposition}
\begin{enumerate}[(i)]
\item For any pretype $A$ we have $FI(\overline{A})=(FI_\bullet(A),FI^\bullet(A))$.
\item If $A, B$ are pretypes and $A\subseteq B$ then $\overline B\subseteq \overline A$.
\item
A pretype $A$ is a type iff $A=\overline{B}$ for some pretype $B$.
\item
The dualizing type $\bbot=\overline{\{1\}}$ is indeed a type. $\Box$
\end{enumerate}
\ep

\subsubsection{Multiplicative operations}
\bd
Let $A$, $B$ be types such that $FI^\bullet(A)\cap FI_\bullet(B)=FI_\bullet(A)\cap FI^\bullet(B)=\emptyset$. The {\it internal hom  type}
$A\multimap B$ is the set
\[
A\multimap B=\{t|~\forall s\in A~ ts\in B\}.
\]
\ed
%\bp
%The pretype $A\multimap B$ is the dual of the pretype $\{ts|~t\in A, s\in \overline{B}\}$. In particular, $A\multimap B$ is a type. $\Box$
%\ep
\bd
Let $A$, $B$ be types such that $FI^\bullet(A)\cap FI^\bullet(B)=FI_\bullet(A)\cap FI_\bullet(B)=\emptyset$. The {\it product pretype} $A\cdot B$, the {\it tensor product type} $A\otimes B$ and the {\it cotensor product type} are given by
$$
\begin{array}{c}
A\cdot B=\{ts|~t\in A, s\in B\},%\\
~A\parr B=\overline{\overline{A}\cdot\overline{B}},%\\
~A\otimes B=\overline{\overline{A\cdot B}},
\end{array}
$$
\ed
Proposition \ref{dual pretypes proposition}, clause (iii) guarantees that tensor and cotensor products of types return types (and not just pretypes).
That the the internal hom  is also a type, follows from the simple proposition below.
\bp\label{cotensor of types}
For all types $A,B$ we have
\begin{enumerate}[(i)]
\item
%$FI(A\cdot B)=
$FI(A\otimes B)=FI(A\parr B)=
((FI^\bullet(A)\cup FI^\bullet(B))\setminus X,
(FI_\bullet(A)\cup FI_\bullet(B))\setminus X)$,
where $X=(FI^\bullet(A)\cap FI_\bullet(B))\cup(FI_\bullet(A)\cap FI^\bullet(B))$, and it is assumed that $A\otimes B$ is defined;
\item  $A\otimes B=\overline{\overline B\parr \overline A}$, $A\parr B=\overline{\overline B\otimes \overline A}$, assumong $A\otimes B$ is defined;
\item $A\multimap B=\overline{A\cdot\overline B}=\overline{A}\parr B$, assuming $A\multimap B$ is defined. $\Box$
\end{enumerate}
\ep
Note that  tensor and cotensor products of types are partial operations. In addition, they are, generally speaking, {\it non-associative}. If we have types $A,B,C$ with $FI^\bullet(A)\cap FI_\bullet(B)\cap FI^\bullet(C)$ nonempty, then, as is seen from the preceding proposition, clause {(i)}, the types $(A\otimes B)\otimes C$ and $A\otimes(B\otimes C)$ (assuming that they are well-defined) even have different boundaries. We will assume the  convention
similar to the one assumed for iterated product of terms (in Definition \ref{iterated multiplication}).
\bd\label{iterated tensor convention}
The expression $X_{1}\otimes\ldots\otimes X_{n}$ denotes the iterated tensor product of types satisfying the condition
$FI^\bullet(X_{i})\cap FI^\bullet(X_{j})=FI_\bullet(X_{i})\cap FI_\bullet(X_{j})=\emptyset$
for all  $i\not=j$. Similarly for cotensor product.
\ed
 In the sequel, when writing tensor and cotensor products of types we will  implicitly assume that the operations are well-defined and,  in the case of iterated products, the convention from Definition \ref{iterated tensor convention} applies. This  makes expressions with iterated products unambiguous.
\bp\label{cotensor properties}
\begin{enumerate}[($i$)]
\item If $A,B,C$ are types and the tensor/cotensor products below are defined then
$$
\begin{array}{c}
A\parr(B\parr C)=\overline{\overline A\cdot\overline B\cdot\overline C},\quad
A\otimes(B\otimes C)=\overline{\overline{A\cdot B\cdot C}}.
\end{array}
$$
\item Tensor and cotensor product of types are associative whenever defined and the convention from Definition \ref{iterated tensor convention} applies.
    \item A term $t$ is in $A_1\parr\ldots\parr A_n$ iff for all
    $$s_{(1)}\in \overline{A_1},\ldots s_{(i-1)}\in \overline{A_{i-1}},
    k_{(1)}\in \overline{A_{i+1}},\ldots k_{(n-i)}\in \overline{A_n}$$ it holds that
    $$ts_{(1)}\ldots s_{(i-1)}k_{(1)}\ldots k_{(n-i)}\in A_i.\quad
    \Box$$
    \end{enumerate}
\ep

\subsubsection{Units and reparameterization}
\bd
Given two indices $x,y$, $x\not=y$, the {\it counit type} $\bot^x_y$  and the  {\it unit type} $\bf {1}^x_y$   are defined by
$$\bot^x_y=\overline{\{\delta^y_x\}},~{\bf 1}^x_y=\overline{\overline{\{\delta^x_y\}}}.$$

Given a type $A$ and indices $x\in FI^\bullet(A)$,
$y\not\in FI^\bullet(A)\cup FI_\bullet(A)$, we define the {\it reparameterization} $A^{[\nicefrac{y}{x}]}$ as $A^{[\nicefrac{y}{x}]}=\{t^{[\nicefrac{y}{x}]}|~t\in A\}$.

Similarly, if
$x\in FI_\bullet(A)$,
$y\not\in FI^\bullet(A)\cup FI_\bullet(A)$, we define
$A_{[\nicefrac{y}{x}]}=\{t_{[\nicefrac{y}{x}]}|~t\in A\}$.
\ed
\bp\label{reparameterization properties}
Given a type $A$ and indices $x\in FI^\bullet(A)$, $y\not\in FI^\bullet(A)\cup FI_\bullet(A)$, we have
\begin{enumerate}[(i)]
\item $A^{[\nicefrac{y}{x}]}=\overline{\overline A_{[\nicefrac{y}{x}]}}$, in particular, $A^{[\nicefrac{y}{x}]}$ is a type;
 \item $A^{[\nicefrac{y}{x}]}={\bf 1}^y_x\otimes A=\bot^y_x\parr A$.
 \end{enumerate}

Similarly,  if $x\in FI_\bullet(A)$, $y\not\in FI^\bullet(A)\cup FI_\bullet(A)$, we have
$A_{[\nicefrac{y}{x}]}=\overline{\overline A^{[\nicefrac{y}{x}]}}$ and
$A_{[\nicefrac{y}{x}]}=A\otimes{1}_y^x=\bot_y^x\parr A$.
\ep
{\bf Proof}
\begin{enumerate}[(i)]
\item Let $t\in \overline{\overline{A}_{[\nicefrac{y}{x}]}}$. Pick $k\in\overline A$ and let
$k'=k_{[\nicefrac{y}{x}]}$.
 Then $k'\in \overline A_{[\nicefrac{y}{x}]}$ and $tk'\in\bbot$.
Let $t'=t^{[\nicefrac{x}{y}]}$.
From term congruence relations (\ref{tensor term relations}) we have that $t'k=t'k'_{[\nicefrac{x}{y}]}=t^{[\nicefrac{x}{y}]}k'_{[\nicefrac{x}{y}]}=tk'$, hence $t'k\in\bbot$. Since $k\in \overline A$ was arbitrary, we get that $t'\in \overline{\overline{A}}=A$. But $t=(t')^{[\nicefrac{y}{x}]}$, so $t\in A^{[\nicefrac{y}{x}]}$. Since $t\in \overline{\overline{A}_{[\nicefrac{y}{x}]}}$ was arbitrary we get $\overline{\overline{A}_{[\nicefrac{y}{x}]}}\subseteq A^{[\nicefrac{y}{x}]}$. The opposite inclusion is proven analogously.
% A part from $x$ and $y$, the terms $k'$ and $t$ have identical free indices with opposite polarities. It follows that $k't$ has only two free indices, namely $x$ and $y$. So $k't$ can be written as $kt=\sigma[w]_x^y$, where $\sigma$ is 1 or purely singular. On the other hand, the term $kt'$ differs from $kt$ only in replacing the upper free occurrence of $y$ with  $x$, i.e.  $kt'=\sigma[w]^x_x=\sigma[w]$. But $kt'\in\bot$, hence
%$\sigma[w]\in\bot$. It follows from Proposition \ref{counit type} that $kt\in\bot^y_x$. Since $k$ was arbitrary, we have shown that for any $k\in \overline A$ the product $kt\in\bot^y_x$.
\item If  $x\in FI^\bullet(t)$ and $y$ does not occur in $t$ then $t^{[\nicefrac{y}{x}]}=t\delta^y_x$. Hence $A^{[\nicefrac{y}{x}]}\subseteq {\bf 1}^y_x \cdot A\subseteq{\bf 1}^y_x\otimes A$.

If $t\in \bot^y_x\parr A$ then by Proposition \ref{cotensor of types}, clause (iii) we have that $\delta_y^xt\in A$. But
$\delta_y^xt=t^{[\nicefrac{x}{y}]}$, hence $t^{[\nicefrac{x}{y}]}\in A$. It follows that $t\in A^{[\nicefrac{y}{x}]}$. Thus $\bot^y_x\parr A\subseteq A^{[\nicefrac{y}{x}]}\subseteq {\bf 1}^y_x\otimes A$.

By the identical reasoning we get $ \overline A\parr\bot^x_y\subseteq \overline A_{[\nicefrac{y}{x}]}\subseteq \overline A\otimes{\bf 1}^x_y$, which yields by duality (i.e. by Proposition \ref{dual pretypes proposition}, clause (ii)) that ${\bf 1}^y_x\otimes A\subseteq A^{[\nicefrac{y}{x}]}\subseteq \bot^y_x\parr A$. $\Box$
\end{enumerate}

\subsection{Interpretation}
\bd\label{valuation}
Given  the  alphabet $\mathit{Lit}_+$  of (positive) literals, a {\it  tensor term model} of the tensor language (intuitionistic or classical)
generated by $\mathit{Lit}_+$ consists of a  semantic type system $(T,\bbot)$
and a valuation $||.||$ mapping every non-unit atomic formula $A$ (i.e. a formula of the form $p^{i_1\ldots i_k}_{j_1\ldots j_k}$, where $p\in \mathit{Lit}_+$)  to a semantic type $||A||$ with the boundary
 $FI(||A||)=(\{i_1\ldots i_k\},\{j_1\ldots j_k\})$ satisfying the conditions
 $$
 \begin{array}{c}
 ||A^{[\nicefrac{i}{j}]}||=||A||^{[\nicefrac{i}{j}]}\mbox{ for }j\in FI^\bullet(A), i\mbox{ fresh},\\
 ||A_{[\nicefrac{i}{j}]}||=||A||_{[\nicefrac{i}{j}]}\mbox{ for }j\in FI_\bullet(A), i\mbox{ fresh}.
 \end{array}
 $$
 The valuation extends to all tensor formulas (intuitionistic or classical) by
  $$
\begin{array}{c}
||{\bf 1}^i_j||={\bf 1}^i_j, ||A\otimes B||=||A||\otimes ||B||,~||A\multimap B||=||A||\multimap||B||,
\\
||\bot^i_j||=\bot^i_j,~
||\overline A||=\overline{||A||},~
||A\parr B||=||A||\parr ||B||.
\end{array}
$$
\ed
It is clear from the discussion in the preceding subsections (in particular, Proposition \ref{cotensor of types}) that the above extension of the valuation is unambiguous.
\bd
A tensor sequent
$A_{(1)},\ldots, A_{(n)}\vdash A$
(i,e, $\vdash \overline{A_{(1)}},\ldots, \overline{A_{(n)}},A$ in the classical format) is {\it valid} in the given model  if one of the two equivalent condition holds:
\begin{enumerate}[(i)]
  \item whenever $t_{(1)}\in{||A_{(1)}||},\ldots,t_{(n)}\in{||A_{(n-1)}||}$ we get $t_{(1)}\cdots t_{(n)}\in||A||$;
  \item $1\in \overline{||A_{(1)}||}\parr\ldots\parr\overline{||A_{(n)}||}\parr||A||$.
\end{enumerate}
\ed
The equivalence of the two above conditions follows from Proposition \ref{cotensor properties}, clause (iii).
%\bd
%Given a set $\Xi$ of tensor sequents (intuitionistic or classical) and a tensor type model  $M$ (of the given language), we say that $M$ is a {\it model of} $\Xi$ if all elements $\Xi$ are valid in $M$.
%\ed
\bt\label{soundness}
%Given a set $\Xi$ of classical, respectively, intuitionistic tensor sequents,
A sequent  derivable in {\bf TTL}, respectively, {\bf ITTL}
 %from elements of $\Xi$
 is valid in any tensor term model.
 % of $\Xi$. In particular, a  derivable (in {\bf TTL} or {\bf ITTL}) sequent is valid in any tensor type model.
\et
{\bf Proof} For the classical case, use induction on derivation using Propositions \ref{cotensor properties} and \ref{reparameterization properties}.
For the intuitionistic case use that the intuitionistic system is a conservative fragment of a classical one. $\Box$

\section{Embedded fragments}
\subsection{Propositional systems of multiplicative linear logic}
\bd
Given a set $\mathit{Lit}_+$ of   {positive literals}
the set ${Fm}$ of {\it classical multiplicative  formulas} over $\mathit{Lit}_+$  is built  according to the grammar
$Fm::=\mathit{Lit}|(Fm\otimes Fm)|(Fm\parr Fm)$, where the set of all literals $\mathit{Lit}=\mathit{Lit}_+\cup\mathit{Lit}_-$ is defined same as for {\bf TTL}. Duality of multiplicative formulas is defined same as for {\bf TTL}.
\ed
There are two natural systems in the multiplicative language, corresponding to (cyclically) ordered and unordered sequents.
\bd
A {\it cyclically ordered sequence} consists of a finite multiset $s$ equipped with a binary relation $\prec$ on $s$ such that there is an enumeration $x_1,\ldots,x_n$ of elements of $s$ satisfying $x_i\prec x_j$ iff $j=i+1(\mod n)$. We say that such an enumeration  {\it induces the cyclic ordering} of $s$.
\ed
\bd
A {\it classical multiplicative}, respectively {\it cyclic multiplicative}, context is a finite  multiset, respectively
cyclically ordered
sequence $\Gamma$ of  classical multiplicative formulas. A {\it classical
multiplicative}, respectively {\it cyclic multiplicative} sequent is an expression of the form $\vdash\Gamma$, where $\Gamma$ is a  classical multiplicative, respectively cyclic multiplicative, context.

The system of {\it classical multiplicative} \cite{Girard_TCS}, respectively {\it cyclic multiplicative linear logic} \cite{Yetter} ({\bf MLL}, respectively {\bf CMLL}), is the set of classical multiplicative, respectively cyclic multiplicative sequents derivable by the rules $({\rm{Id}})$, $({\rm{Cut}})$,  $(\parr)$ and $(\otimes)$ of {\bf TTL}.
%The system of  {\it cyclic multiplicative logic} ({\bf CMLL}) is the set of classical multiplicative sequents derivable by the rules $({\rm{Id}})$, $({\rm{Cut}})$, $(\parr)$ and $(\otimes)$ of {bf TTL} and the following {\it Rotation} rule
%\[
%\cfrac{\vdash\Gamma,\Delta}{\vdash\Delta,\Gamma}({\rm{Rot}}).
%\]
     \ed
{\it Lambek calculus} (more precisely,  Lambek calculus allowing {\it empty antecedents} \cite{Lambek_empty}) is just the intuitionistic fragment of {\bf CMLL}.
\bd
Given a set $\mathit{Lit}$ of   {literals}
the set ${Tp}$ of {\it Lambek types} over $\mathit{Lit}$  is built  according to the grammar
$Tp::=\mathit{Lit}|(Tp\bullet Tp)|(Tp\backslash TP)|(Tp/Tp)$.
{\it Lambek context} is a finite sequence of Lambek types, {\it Lambek sequent} is an expression of the form $\Gamma\vdash A$, where $\Gamma$ is a Lambek context and $A$ is a Lambek type.

The translation of Lambek types  to {\bf CMLL} is given by
$$A\bullet B=A\otimes B,~ A\backslash B=\overline A\parr B,~ B/A=B\parr\overline A.$$
 The Lambek sequent $A_1,\ldots,A_n\vdash B$ translates to {\bf CMLL} as $\vdash \overline{A_1},\ldots, \overline{A_n},B$.
 {\it Lambek calculus} ({\bf LC}) is the set of Lambek sequents derivable in {\bf CMLL} under the above translation.
\ed
The rules of {\it Lambek calculus} ({\bf LC}) that can be found, for example, in \cite{Lambek_empty}.

\subsection{{\bf TTL} and noncommutative logic}
\subsubsection{Embedding cyclic logic in {\bf TTL}}
Let the set $\mathit{Lit}_+$  of positive literals be given. We equip each element $p\in\mathit{Lit}_+$ with valency $(1,1)$ and define a translation from the  classical multiplicative language over $\mathit{Lit}_+$ to the classical tensor language over $\mathit{Lit}_+$, which will induce conservative embedding of {\bf CMLL} and {\bf LC}.
\bd
In the setting as above, the {\it cyclic translation} $||F||^i_j$ parameterized by $i,j\in\mathit{Ind}$, where $i\not=j$, from a classical multiplicative formula $F$ to the tensor language is defined up to renaming bound indices by induction as
\be\label{cyclic translation}
||p||^i_j=p^i_j\mbox{ for }p\in\mathit{Lit},~||A\otimes B||^i_j=||A||^i_k\otimes||B||^k_j,~||A\parr B||^i_j=||A||^i_k\parr||B||^k_j,
\ee
where the representatives of $||A||^i_k$ and $||B||^k_j$ are chosen so that they have no common indices except $k$.
\ed
\bp
If a classical multiplicative sequent $\vdash A_{1},\ldots,A_n$, where $n>1$, is {\bf CMLL} derivable
and $i_1,\ldots,i_n$ are indices such that $\Gamma=||A_1||^{i_0}_{i_1},||A_2||^{i_1}_{i_2},\ldots,||A_n||^{i_{n-1}}_{i_0}$ is a well-formed tensor context
then the sequent
$\vdash\Gamma$ is derivable in {\bf TTL}.
If $\vdash A$ is {\bf CMLL} derivable then the sequent $\vdash\bot^i_j,||A||^j_i$ is derivable in {\bf TTL}.
\ep
{\bf Proof} Induction on a cut-free derivation in {\bf CMLL}. All steps are shown in Figure \ref{Translating  CMLL derivations}. $\Box$
\begin{figure}
  \centering
  %%%%%%%%%%%%%%%%%%%%%%%%%%%%%%%%%%%%%%%%
\subfloat[{Translating the $(\otimes)$ rule}
%\label{side cut figure2}
]
{
  $
  \begin{array}{c}
      \def\ScoreOverhang{.1pt}
\AxiomC{$\vdash A_1,\ldots,A_m$~$m>1$}
\AxiomC{$\vdash B_1,\ldots,B_n$~$n>1$}
\RightLabel{($\otimes$)}
\insertBetweenHyps{\hskip 1pt}
\BinaryInfC{$\vdash A_1,\ldots,A_{m-1},A_m\otimes B_1,B_2,\ldots,B_n$}
\DisplayProof
\Longrightarrow
\\[.3cm]
  \def\ScoreOverhang{.1pt}
\AxiomC{$\vdash||A_1||^{i_0}_{i_1},\ldots,||A_m||^{i_{m-1}}_{i_0}$}
\AxiomC{$\vdash||B_1||^{i_0}_{i_{n}},\ldots,||B_n||^{i_{n+m-1}}_{i_{0}}$}
\RightLabel{($\otimes$)}
\insertBetweenHyps{\hskip 1pt}
\BinaryInfC{$\vdash||A_1||^{i_0}_{i_1},\ldots,||A_{m-1}||^{i_{m-2}}_{i_{m-1}},||A_m||^{i_{m-1}}_{i_0}\otimes||B_1||^{i_0}_{i_{n}},
||B_2||^{i_{n}}_{i_{n+1}}\ldots,||B_n||^{i_{n+m-1}}_{i_{0}}$}
\DisplayProof\\[.7cm]
  \def\ScoreOverhang{.1pt}
\AxiomC{$\vdash A$}
\AxiomC{$\vdash B_1,\ldots,B_n$~$n>1$}
\RightLabel{($\otimes$)}
\insertBetweenHyps{\hskip 5pt}
\BinaryInfC{$\vdash A\otimes B_1,\ldots,A_{n-1},A_n$}
\DisplayProof
\Longrightarrow
%\\[.3cm]
  \def\ScoreOverhang{.1pt}
\AxiomC{$\vdash||A||^{i_0}_j,\bot^j_{i_0}$}
\AxiomC{$\vdash||B_1||^{j}_{i_1},\ldots,||B_n||^{i_{n-1}}_{j}$}
\RightLabel{($\otimes$)}
\insertBetweenHyps{\hskip 5pt}
\BinaryInfC{$\vdash||A||^{i_0}_j\otimes||B_1||^{j}_{i_1},||B_2||^{i_1}_{i_2},\ldots,||B_n||^{i_{n-1}}_{j},\bot^{j}_{i_{0}}$}
\RightLabel{($\bot\leftarrow$)}
\UnaryInfC{$\vdash||A||^{i_0}_j\otimes||B_1||^{j}_{i_1},||B_2||^{i_1}_{i_2},\ldots,||B_n||^{i_{n-1}}_{i_0}$}
\DisplayProof\\[.7cm]
 \def\ScoreOverhang{.1pt}
\AxiomC{$\vdash A$}
\AxiomC{$\vdash B$}
\RightLabel{($\otimes$)}
%\insertBetweenHyps{\hskip 5pt}
\BinaryInfC{$\vdash A\otimes B$}
\DisplayProof
\Longrightarrow
  \def\ScoreOverhang{.1pt}
\AxiomC{$\vdash||A||^{j}_k,\bot^k_{j}$}
\AxiomC{$\vdash||B||^{k}_{i},\bot^i_k$}
\RightLabel{($\otimes$)}
\BinaryInfC{$\vdash||A||^{j}_k\otimes||B||^{k}_{i},\bot^{k}_{j},\bot^i_k$}
\RightLabel{($\bot\leftarrow$)}
\UnaryInfC{$\vdash||A||^{j}_k\otimes||B||^{k}_{i},\bot_{j}^{i}$}
\DisplayProof
\end{array}
$
}\\
 %%%%%%%%%%%%%%%%%%%%%%%%%%%%%%%%%%%%%%%%
\subfloat[{Translating the $(\parr)$ rule}
%\label{side cut figure2}
]
{
  $
  \begin{array}{c}
      \def\ScoreOverhang{.1pt}
\AxiomC{$\vdash A_1,\ldots,A_n$~$n>2$}
\RightLabel{($\parr$)}
\UnaryInfC{$\vdash A_1,\ldots,A_{n-2},A_{n-1}\parr A_n$}
\DisplayProof
\Longrightarrow
\\[.3cm]
  \def\ScoreOverhang{.1pt}
\AxiomC{$\vdash||A_1||^{i_0}_{i_1},\ldots,||A_{n-1}||^{i_{n-2}}_{i_{n-1}},||A_n||^{i_{n-1}}_{i_0}$}
\RightLabel{($\parr$)}
\UnaryInfC{$\vdash||A_1||^{i_0}_{i_1},\ldots,||A_{n-2}||^{i_{n-3}}_{i_{n-2}},||A_{n-1}||^{i_{n-2}}_{i_{n-1}}\parr||A_n||^{i_{n-1}}_{i_0}$}
\DisplayProof\\[.7cm]
  \def\ScoreOverhang{.1pt}
\AxiomC{$\vdash A,B$}
\RightLabel{($\parr$)}
\UnaryInfC{$\vdash A\parr B$}
\DisplayProof
\Longrightarrow
  \def\ScoreOverhang{.1pt}
\AxiomC{$\vdash||A||^{i}_k,||B||^k_i$}
\RightLabel{($\bot^\to$)}
\UnaryInfC{$\vdash\bot^i_j,||A||^{j}_k,||B||^{k}_{i}$}
\RightLabel{($\parr$)}
\UnaryInfC{$\vdash\bot^i_j,||A||^{j}_k\parr||B||^{k}_{i}$}
\DisplayProof
 \end{array}
$
}
  \caption{Translating {\bf CMLL} derivations}
 \label{Translating  CMLL derivations}
\end{figure}
\bd
A tensor formula is {\it cyclic} if it is in the image of translation (\ref{cyclic translation}). The {\it cyclic fragment} of {\bf TTL} consists of sequents all whose formulas are cyclic. The {\it extended cyclic fragment} of {\bf TTL} consists of sequents all whose formulas are cyclic or counits $\bot^i_j$.
\ed
\bd
For a tensor sequent $\vdash\Gamma$ in the extended cyclic fragment we define the {\it neighboring relation}
$\prec$
 between formulas of $\Gamma$ by $A\prec B$ iff $FI_\bullet(A)=FI^\bullet(B)$.
\ed
\bl
For any {\bf TTL} derivable sequent in the extended cyclic fragment the neighboring relation is a cyclic ordering.
\el
{\bf Proof} Induction on a cut-free derivation. As an illustration consider the induction step corresponding to the $(\otimes)$ rule.

Let the sequent $\vdash\Theta$ of the form $\vdash\Gamma,A\otimes B,\Delta$ in the extended cyclic fragment be obtained from {\bf TTL} derivable sequents
$\vdash\Theta_1$, $\vdash\Theta_2$ of the form $\vdash \Gamma,A$, $\vdash B,\Delta$ respectively by the $(\otimes)$ rule.

 By the induction hypothesis we have  enumerations of $\Theta_1$ and $\Theta_2$, respectively
$A_{(1)},\ldots,A_{(n)},A$ and $B, B_{(1)},\ldots,B_{(m)}$,  that induce the respective neighboring relations. In particular, we have
$FI_\bullet(A_{(i)})=FI^\bullet( A_{(i+1)})$  for $i<n$,
$FI_\bullet(A_{(n)})=FI^\bullet( A)$, $FI_\bullet(A)=FI^\bullet( A_{(1)})$
and similarly
$FI_\bullet(B_{(i)})=FI^\bullet( B_{(i+1)})$ for $i<m$,
$FI_\bullet(B_{(m)})=FI^\bullet(B)$, $FI_\bullet(B)=FI^\bullet( B_{(1)})$.

Now $A\otimes B$ must be a cyclic formula, which means $FI_\bullet(A)=FI^\bullet(B)$. Then we compute
$FI_\bullet(A\otimes B)=FI_\bullet(B)=FI^\bullet(B_{(1)})$ and
 $FI^\bullet(A\otimes B)=FI^\bullet(A)=FI_\bullet(A_{(n)})$.
 Finally $FI_\bullet(B_m)=FI^\bullet(B)=FI_\bullet(A)=FI^\bullet(A_{(1)})$.

 It follows that the neighboring relation on $\Gamma,A\otimes B,\Delta$ coincides with the cyclic ordering induced by the enumeration
 $A_{(1)},\ldots,A_{(n)}, A\otimes B,B_{(1)},\ldots,B_{(m)}$. $\Box$
\bp
Let $\vdash \Gamma$ be a {\bf TTL} derivable sequent in the extended cyclic fragment, and let
$A_{(1)},\ldots A_{(n)}$ be an enumeration of formulas in $\Gamma$ inducing the neighboring relation.
 Let the cyclic sequent $\vdash\Gamma'$ be obtained from
the expression $\vdash A_{(1)},\ldots A_{(n)}$  by erasing all indices and all counits. Then $\vdash\Gamma'$ is derivable in {\bf CMLL}.
\ep
{\bf Proof} Induction on a cut-free derivation of $\vdash\Gamma$. $\Box$
\bc\label{cyclic translation conservative}
The cyclic translation is conservative.
The {\bf CMLL} sequent $\vdash A_{1},\ldots,A_n$, where $n>1$,  is derivable (in {\bf CMLL}) iff for some (hence any) array of indices
 $i_1,\ldots,i_n$  such that $\Gamma=||A_1||^{i_1}_{i_2},||A_2||^{i_2}_{i_3},\ldots,||A_n||^{i_n}_{i_1}$ is a well-formed tensor context
the sequent
$\vdash\Gamma$ is derivable in {\bf TTL}. The {\bf CMLL} sequent $\vdash A$ is derivable in {\bf CMLL} of the sequent $\vdash\bot^i_j,||A||^j_i$ is derivable in {\bf TTL}. $\Box$
\ec

\subsubsection{Embedding Lambek calculus}
Given that {\bf LC} is a conservative fragment of {\bf CMLL} and {\bf ITTL} is a conservative fragment of {\bf TTL}, the embedding of {\bf CMLL} into {\bf TTL} induces the embedding of {\bf LC} into {\bf ITTL}.
As in the case of {\bf CMLL}, we equip each literal $p\in\mathit{Lit}$ with valency $(1,1)$ and define a translation from the language of {\bf LC} over $\mathit{Lit}$ to the intuitionistic tensor language over $\mathit{Lit}$.
\bd
In the setting as above, the {\it translation} $||F||^i_j$ parameterized by $i,j\in\mathit{Ind}$, where $i\not=j$, from a Lambek type $F$ to the intuitionistic tensor language is defined up to renaming bound indices by induction as
$$
\begin{array}{rl}
||p||^i_j=p^i_j,&||A\bullet B||^i_j=||A||^i_k\otimes||B||^k_j,\\
||A\backslash B||^i_j=||A||^k_i\multimap||B||^k_j,&
||B/ A||^i_j=||A||_k^j\multimap||B||_k^i,
\end{array}
$$
where the representatives in every pair $(||A||^i_k,||B||^k_j)$,
$(||A||^k_i,||B||^k_j,)$,
$(||A||_k^j,||B||_k^i)$
 are chosen so that they have no common indices except $k$.
\ed
Conservativity of {\bf LC} in {\bf CMLL} and {\bf ITTL} in {\bf TTL} (together with the provable equivalence
$\overline{||A||_k^j}\parr||B||_k^i\cong ||B||_k^i\parr\overline{||A||_k^j}
)$
immediately yield the following corollary from Corollary \ref{cyclic translation conservative}
\bc
The {\bf LC} sequent $A_{1},\ldots,A_n\vdash A$, where $n>0$,  is derivable (in {\bf LC}) iff for some (hence any) array of indices
 $i_1,\ldots,i_n$  such that $\Gamma=||A_1||^{i_0}_{i_1},\ldots,||A_n||^{i_{n-1}}_{i_n}$ is a well-formed intuitionistic tensor context
the sequent
$\Gamma\vdash||A||^{i_0}_{i_n}$ is derivable in {\bf ITTL}.
The sequent $\vdash A$ is derivable in {\bf LC} iff ${1}^i_j\vdash||A||^i_j$ is derivable in {\bf ITTL}. $\Box$
\ec
\begin{Rem}\label{Rem3}
Comparing the above translation from {\bf LC} to {\bf ITTL} with formulas (\ref{Lambek2MILL1}) giving the translation from {\bf LC} to {\bf MILL1} suggests that {\bf ITTL} could itself be translated to {\bf MILL1}, with indices becoming first order variables and bound indices corresponding to variables bound by a quantifier, according to the informal scheme
\[A^i\otimes B_i\mapsto \exists i(A(i)\otimes B(i)),~A^i\multimap B^i\mapsto \forall i(A(i)\multimap B(i)).\]
 As for the propositional constant ${\bf 1}$, the {\bf ITTL} formula ${\bf 1}^i_j$ behaves much like the equality $(i=j)$, so one might attempt a translation into the first order language with equality.

 We note however that such a translation could not be conservative. Indeed, the sequent $\vdash (A^i_j\multimap A^i_k)\otimes(A^l_k\multimap A^l_j)$, which is underivable in {\bf ITTL} (see Remark \ref{Rem2}), should translate as
$\vdash \exists k\exists j(\forall i(A(i,j)\multimap A(i,k))\otimes\forall l(A(l,k)\multimap A(l,j)))$, and the latter is {\bf MILL1} derivable. It seems that in order to have a direct relation with {\bf MILL1} we need to consider tensor terms with congruence (\ref{empty loop}) added to the definition (see Remarks \ref{Rem1}, \ref{Rem2}).
\end{Rem}

%%%%%STOP HERE
\subsection{TTL and MLL}
\subsubsection{{\bf MLL}-like fragment}
\bd
A tensor formula is {\it {\bf MLL}-like} if it does not involve units ${\bf 1}^i_j$ or counits $\bot^i_j$ and has no bound indices. The {\it {\bf MLL}-like fragment} of {\bf TTL} consists of sequents all whose formulas are {\bf MLL}-like.
\ed
\bp\label{derivations in MLL-like}
 If an  {\bf MLL}-like sequent $\vdash\Gamma$ is derivable in  {\bf TTL} then it is derivable not using the $(\bot_\leftarrow)$ and the $(\bot^\to)$ rules.
\ep
{\bf Proof} By Corollary \ref{essential reductions}, if $\vdash\Gamma$ has a derivation in {\bf TTL} then it has a cut-free derivation where all applications of the $(\bot_\leftarrow)$ rule are essential. But a cut-free derivation of $\vdash\Gamma$ involves only subformulas of $\Gamma$ and, possibly, counits, and these have no bound indices. Without bound indices we cannot have essential applications of the $(\bot_\leftarrow)$ rule, hence there are no $(\bot_\leftarrow)$ rules at all. Now, $\vdash\Gamma$ does not contain counits, and since the derivation is cut-free and does not use the $(\bot_\leftarrow)$ rule, it does not use the $(\bot^\to)$ rule, which introduces counits, either. $\Box$

We see that the {\bf MLL}-like fragment uses only  rules of {\bf MLL}, which suggests that the above fragment and the system of {\bf MLL} should be translatable to each other. And indeed there exists a translation from {\bf MLL} to the {\bf MLL}-like fragment, however the input of the translation is not a provable {\bf MLL} sequent, but a pair: an {\bf MLL} sequent together with its derivation. ``Morally'' the situation is as follows. Given an {\bf MLL} sequent $\vdash\Gamma$, we try to decorate it with indices to obtain a {\bf TTL} derivable sequent from the {\bf MLL}-like fragment. Different successful decorations of $\vdash\Gamma$ correspond to its different derivations in {\bf MLL} (more precisely, to different proof-nets).

\subsubsection{Embedding MLL derivations}
In this section we will  follow notational conventions of Section \ref{cut-elimination section}, where we used capital Latin letters for sequences of indices. Also, we will write $\bot^I_J$ for the context $\bot^{i_1}_{j_1},\ldots,\bot^{i_n}_{j_n}$, where $I=i_1,\ldots,i_n$,
$J=j_1,\ldots,j_n$.

We assume that we are given the alphabet $\mathit{Lit}_+$ of positive literals generating a classical multiplicative language, and that elements of $\mathit{Lit}_+$ are, furthermore, assigned valencies so that they generate a classical tensor language as well.
\bd
For an {\bf MLL} sequent $\Pi$ of the form $\vdash A_1\ldots A_n$ together with a derivation $\pi$ of $\Pi$  the {\it {\bf TTL} representation} of the pair $(\Pi,\pi)$ is a {\bf TTL} derivable tensor sequent $||(\Pi,\pi)||$ of the form $\vdash ||A_1||^\pi,\ldots,||A_n||^\pi$, where each tensor formula $||A_i||^\pi$ is $A_i$ decorated with some indices,  defined up to similarity of tensor sequents by induction on $\pi$ as follows.
\begin{enumerate}[(i)]
\item If $\pi$ is an axiom and $\Pi$ has the form $\vdash \overline p,p$, where $v(p)=(m,n)$, choose two sequences of pairwise distinct indices
$I=i_1,\ldots,i_m$, $J=j_1,\ldots,j_m$, $I\cap J=\emptyset$ and put
$||p||^\pi=p^{I}_{J}$, $||\overline p||^\pi=\overline{||p||^\pi}$. The obtained sequent  is an instance of the $({\rm{Id}})$ axiom of {\bf TTL}.
%\item If
%$\pi$ with the conclusion $\vdash \Gamma,B,A,\Delta$ is obtained from a derivation $\pi'$ with the conclusion $\vdash\Gamma,A,B,\Delta$ by the $({\rm{Ex}})$ rule,  then the sequent
%$||\pi||$ is the one obtained from $||\pi'||$ by permuting the formulas $||A||^{\pi'}$ and $||B||^{\pi'}$ (and derivable from $||\pi'||$ by the $({\rm{Ex}})$ rule.
\item If $\Pi$ has the form $\vdash \Gamma,A\parr B$ and $\pi$ is obtained from a derivation $\pi'$ of the sequent $\Pi'$ having the form $\vdash\Gamma,A,B$ by the $(\parr)$ rule, put $||F||^\pi=||F||^{\pi'}$ for $F\in\Gamma$ and $||A\parr B||^\pi=||A||^{\pi'}\parr||B||^{\pi'}$. The sequent $||(\Pi,\pi)||$ is derivable from $||(\Pi',\pi')||$ by the $(\parr)$ rule.
\item If $\Pi$ has the form $\vdash \Gamma,A\otimes B,\Delta$ and $\pi$ is obtained from  derivations $\pi_1$, $\pi_2$ with the conclusions
$\Pi_1$, $\Pi_2$ of the form  $\vdash\Gamma,A$, $\vdash B,\Delta$ respectively by the $(\otimes)$ rule,
    choose representatives of $||(\Pi_1\pi_1)||$ and $||(\Pi_2,\pi_2)||$ that have no indices in common and put
     $||F||^\pi=||F||^{\pi_1}$ for $F\in\Gamma$,
    $||F||^\pi=||F||^{\pi_2}$ for $F\in\Delta$
    and $||A\otimes B||^\pi=||A||^{\pi_1}\otimes||B||^{\pi_2}$. The sequent $||(\Pi,\pi)||$ is derivable from $||(\Pi_1,\pi_1)||$ and $||(\pi_2,\pi_2)||$ by the $(\otimes)$ rule.
\item If $\Pi$ has the form $\vdash \Gamma,\Delta$ and $\pi$ is obtained from  derivations $\pi_1$, $\pi_2$ with the conclusions
$\Pi_1$, $\Pi_2$ of the form $\vdash\Gamma,A$, $\vdash \overline A,\Delta$ respectively by the Cut rule,
    choose representatives of $||(\Pi_1,\pi_1)||$ and $||(\Pi_2,\pi_2)||$ that have no indices in common. Let
    $FI(||A||^{\pi_1})=(I,J)$ and $FI(||\overline A||^{\pi_2})=(K,L)$.
         Put $||F||^\pi=(||F||^{\pi_1})_{[\nicefrac{L}{I}]}$ for $F\in\Gamma$,
    $||F||^\pi=(||F||^{\pi_2})_{[\nicefrac{J}{K}]}$ for $F\in\Delta$.
    The sequent $||\pi||$ is derivable from $||\pi_1||$ and $||\pi_2||$ by
    a combination of the rules $(\bot^\to)$, Cut and $(\bot_\leftarrow)$ as shown in Figure \ref{translating cut}.
\end{enumerate}
\ed
\begin{figure}
  \centering
  $
  \begin{array}{c}
  FI(||A||^{\pi_1})=(I,J),~FI(||\overline A||^{\pi_2})=(K,L)\\
\def\ScoreOverhang{.1pt}
\AxiomC{$\vdash||\Gamma||^{\pi_1},\fCenter||A||^{\pi_1}$}
\RightLabel{($\bot^\to$)}
\def\extraVskip{8pt}
\UnaryInfC{$ \cdots$}
\RightLabel{($\bot^\to$)}
\def\extraVskip{3pt}
\UnaryInfC{$ \vdash||\Gamma||^{\pi_1},\bot^I_L,(||A||^{\pi_1})^{[\nicefrac{L}{I}]}$}
\def\extraVskip{2pt}
\AxiomC{$\vdash||\overline A||^{\pi_2},||\Delta||^{\pi_2}$}
\RightLabel{($\bot^\to$)}
\UnaryInfC{$\cdots$}
\RightLabel{($\bot^\to$)}
\UnaryInfC{$ \vdash\bot^K_J,(||\overline A||^{\pi_2})^{[\nicefrac{J}{K}]},||\Delta||^{\pi_2}$}
\def\extraVskip{2pt}
\RightLabel{(${\rm{Ex}}$)}
\UnaryInfC{$ \vdash(||\overline A||^{\pi_2})^{[\nicefrac{J}{K}]},\bot^K_J,||\Delta||^{\pi_2}$}
\RightLabel{(${\rm{Cut}}$)}
%\insertBetweenHyps{\hskip 1pt}
\BinaryInfC{$ \vdash ||\Gamma||^{\pi_1},\bot^I_L,\bot^K_J,||\Delta||^{\pi_2}$}
\RightLabel{($\bot_\leftarrow$)}
\def\extraVskip{4pt}
\UnaryInfC{$ \cdots$}
\RightLabel{($\bot_\leftarrow$)}
\def\extraVskip{2pt}
\UnaryInfC{$ \vdash (||\Gamma||^{\pi_1})_{[\nicefrac{L}{I}]},(||\Delta||^{\pi_2})_{[\nicefrac{J}{K}]}$}
\DisplayProof
\end{array}
$
  \caption{Translating {\bf MLL} cut to {\bf TTL}}
  \label{translating cut}
\end{figure}
\bp
{\bf TTL} translation of  {\bf MLL} is an invariant of cut-elimination. If an {\bf MLL} derivation $\pi'$ of a sequent $\Pi$ is obtained from the derivation $\pi$ by cut-elimination then $||(\Pi,\pi)||$ and $||(\Pi,\pi')||$ coincide.
\ep
{\bf Proof} Induction on the cut-elimination algorithm for {\bf MLL}. $\Box$
\bt
Every {\bf TTL} derivable sequent $\vdash\Gamma$ from the {\bf MLL}-like fragment  is the translation  of an {\bf MLL} sequent with a chosen derivation.
\et
{\bf Proof} Induction on a cut-free derivation of $\vdash\Gamma$ not using the $(\bot_\leftarrow)$ and $(\bot^\to)$ rules. $\Box$

It is quite obvious that we can restrict the above construction to the intuitionistic fragment of {\bf MLL} and obtain a translation to {\bf ITTL}.
Now we note that, in the intuitionistic setting, the inputs of our translation, i.e. sequents equipped with chosen derivations, canonically correspond via the Curry-Howard isomorphism to $\lambda$-calculus typing judgements, and it is precisely $\lambda$-calculus typing judgements that play the main role in the ``commutative'' categorial grammars, i.e. {\bf ACG}. On one hand, this suggests that, as far as categorial grammars are concerned, the input of the translation is chosen correctly. On the other hand this leads to the question of representing {\bf ACG} in {\bf ITTL}. We leave this question to future research.

\section{Grammars}%\label{lexcions section}
%In this section we discuss grammars based on tensor term logic as the typing system. These may be intuitionistic or classical depending on the chosen logic, many definitions and statements are identical modulo this difference. We will sometimes omit explicit reference to the particular typing system, implying that the definition or statement can be read as two versions, respectively classical and intuitionistic.
%
\subsection{Typing judgements and grammars}
\bd  A classical, respectively intuitionistic {\it typing judgement}  is an expression of the form $t\rhd A$, where $t$ is a term and $A$ is a classical, respectively intuitionistic formula with the same boundary: $FI(t)=FI(A)$.
The typing judgement is  {\it regular}, respectively {\it lexical}, if the term $t$ is regular, respectively lexical.
 A typing judgement  of the form $\delta^i_j\rhd {\bf 1}^i_j$ is a {\it typing axiom}. The set of all typing axioms is denoted as $\mathit{TpAx}$.
\ed
The intended  semantics of the typing judgement $t\rhd A$ is that $t\in||A||$. With such a semantics, typing axioms are always valid, and this explains the title ``axiom''.

Typing judgements can be understood as a special kind of sequents and grammars as systems of nonlogical axioms (this will be discussed shortly below). In particular, we can define similarity and  $\alpha$-equivalence of typing judgements.
\bd
{\it Similarity of typing judgements} is generated by the following relations:
\[
\begin{array}{c}
j\in FI^\bullet(t), i\mbox{ fresh}\Rightarrow (t^{[\nicefrac{i}{j}]}\rhd A^{[\nicefrac{i}{j}]})\sim (t\rhd A)\\
j\in FI_\bullet(t), i\mbox{ fresh}\Rightarrow (t_{[\nicefrac{i}{j}]}\rhd A_{[\nicefrac{i}{j}]})\sim (t\rhd A).
\end{array}
\]
{\it $\alpha$-Equivalence of typing judgements} is generated by the following relations:
\[
\begin{array}{c}
(t\rhd A)\sim (t\rhd B)\Rightarrow (t\rhd A)\equiv_\alpha (t\rhd B),\\
A\equiv_\alpha B\Rightarrow (t\rhd A)\equiv_\alpha (t\rhd B).
\end{array}
\]
\ed
\bd
  A  classical, respectively intuitionistic,
  {\it lexicon} $\mathit{Lex}$,    is a set of classical, respectively intuitionistic,
   lexical typing judgements.  The lexicon is {\it regular} if all its elements are regular typing judgements. The {\it closure}  $\mathit{Lex}_{cl}$ of the lexicon $\mathit{Lex}$ is the set $\mathit{Lex}_{cl}=\mathit{Lex}_\sim\cup\mathit{TpAx}$, where $\mathit{Lex}_\sim$ is the closure of $\mathit{Lex}$ under similarity.
  %closed under $\alpha$-equivalence and  containing finitely many $\alpha$-equivalence classes,

A classical, respectively intuitionistic {\it tensor grammar} $G$ is a pair $G=(\mathit{Lex},S)$,    where
$\mathit{Lex}$ is a  finite regular classical, respectively intuitionistic lexicon
and
    $S$, the {\it sentence type symbol},  is a positive literal of valency $(1,1)$.

A  regular typing judgement $t\rhd A$ is {\it derivable in} $G$, notation $t\rhd_GA$, if there exist elements $t_{(1)}\rhd A_{(1)},\ldots,t_{(n)}\rhd A_{(n)}$ of $\mathit{Lex}_{cl}$  such that $t=t_{(1)}\cdots t_{(n)}$ and
the sequent $A_{(1)},\ldots A_{(n)}\vdash A$
(i.e., in the classical case,
   $\vdash\overline{A_{(1)}},\ldots \overline{A_{(n)}}, A$)
is derivable  in {\bf TTL}, respectively {\bf ITTL}. (The number $n$ may equal zero, which means that $t=1$ and the sequent $\vdash A$ is derivable in {\bf TTL},
  respectively {\bf ITTL}.
)

  The {\it language $L(G)$ generated by }$G$ or, simply, the {\it language $L(G)$ of} $G$ is the set $L(G)=\{w\in T^*|~[w]^i_j\rhd_GS^i_j\}$.
  \ed
 %\bp
%The set of typing judgements derivable in a given tensor grammar $G$ (intuitionistic or classical) is closed under $\alpha$-equivalence.
%\ep
%{\bf Proof}
%We will discuss the classical case.
%
%The proof is by induction on the definition of $\alpha$-equivalence. We will discuss one step.
%Let $\mathit{Lex}$ be the lexicon of $G$ and assume, for example, that $t\rhd_GA$, $j\in FI^\bullet(A)$, $i$ is fresh and $t'=t^{[\nicefrac{i}{j}]}$, $A'=A^{[\nicefrac{i}{j}]}$. We need to show that $t'\rhd_GA'$.
%
%By definition, we have axioms $t_{(1)}\rhd A_{(1)},\ldots, t_{(n)}\rhd A_{(n)}$ in $\mathit{Lex}\cup\mathit{TpAx}$ such that $t=t_{(1)}\ldots t_{(n)}$ and the sequent $\vdash\overline{A_{(1)}},\ldots,\overline{A_{(n)}},A$ is derivable in {\bf TTL}.
%
%By the $(\bot^\to)$ rule we get that the sequent
%$\vdash\overline{A_{(1)}},\ldots,\overline{A_{(n)}},\bot_i^j,A'$ is derivable in {\bf TTL} as well. Also we have $t'=t\delta^i_j=t_{(1)}\ldots t_{(n)}\delta^i_j$. Put $t_{(n+1)}=\delta^i_j$, $A_{(n+1)}={1}^i_j$. Then $t_{(n+1)}\rhd A_{(n+1)}$ is a typing axiom and we get the list
%$t_{(1)}\rhd A_{(1)},\ldots, t_{(n+1)}\rhd A_{(n+1)}$ in $\mathit{Lex}\cup\mathit{TpAx}$ such that $t'=t_{(1)}\ldots t_{(n+1)}$ and the sequent $\vdash\overline{A_{(1)}},\ldots,\overline{A_{(n+1)}},A'$ is derivable in {\bf TTL}. Which means exactly that $t'\rhd_GA'$. $\Box$
Semantically, the typing judgements derivable in a tensor grammar are those that are valid whenever the elements of the lexicon are valid. This applies, in particular, to elements of the lexicon closure.
%\bp
%Let $G$ be a tensor grammar and assume that $t\rhd_GA$. Let $G'$ be the grammar obtained from $G$ by adding the typing judgement $t\rhd A$ to the lexicon. Then $G'$ has the same set of derivable typing judgements as $G$. $\Box$
%\ep

Prior to giving an example of a tensor grammar we will discuss a convenient representation of those as systems of sequents, {\it non-logical axioms}, in the extended tensor formula language where terminal symbols are treated as additional literals.

\subsection{Typing judgements as sequents}\label{extended language}
\bd
Given a  set $\mathit{Lit}_+$ of  positive literals with assigned valencies and a terminal alphabet $T$, the {\it extended} tensor language is the tensor language (respectively intuitionistic or classical)
generated by $\mathit{Lit}_+\cup T$, where every element of $T$ is assigned valency $(1,1)$. The subset of the extended language not using terminal symbols (i.e. generated by $\mathit{Lit}_+$) is the {\it base} language.
\ed
When writing formulas and sequents in the extended language, it will be convenient to distinguish notationally terminal symbols and formulas of the base language. When using an element $a\in T$ as a literal in the extended language we will write it in boldface as ${\bf a}$.

 The extended language allows representing terms over $T$ as formulas and  typing judgements as sequents. Typically, ``atomic'' terms of the form $[a]^i_j$ can be seen as atomic formulas of the form ${\bf a}^i_j$, the Kronecker delta $\delta^i_j$ corresponds to the unit ${\bf 1}^i_j$ and the product of terms corresponds to the tensor product of formulas.
 For the  proof-theoretic analysis,  however, it is more convenient to represent terms not as formulas, but as {\it tensor contexts} in the extended language.
 \bd
 The map $\mathit{ct}$ from well-formed pseudo-terms to {\it classical} tensor contexts in the extended language is {\it partially} defined by the following induction:
 \[
\mathit{ct}(1)=\emptyset, \mathit{ct}([a]^i_j)={\bf \overline a}_i^j\mbox{ for }a\in T,~\mathit{ct}(\delta^i_j)=\bot^j_i,
 ~\mathit{ct}(t\cdot s)=\mathit{ct}(t),\mathit{ct}(s).
 \]
 \ed
  \bp\label{surjectivity of tm}
 For every term $t$ over $T$ there exists a well-formed pseudo-term $\tau$ representing $t$ such that $\mathit{ct}(\tau)$ is defined.
 \ep
 {\bf Proof} This follows from the fact that every term can be represented as  the product of elementary regular terms of the form $[a]^i_j$ or $\delta^i_j$, $i\not=j$.

 For the elementary regular term $t=[a_{1}\ldots a_{n}]^i_j$, where
  $a_{1},\ldots, a_{n}\in T$
  $n>0$, we have
 $t=[a_{1}]^i_{i_1}\cdots[a_{n}]^{i_{n-1}}_{j}$, where
$i,i_1,\ldots,i_{n-1},j$ are pairwise distinct.

 For the elementary singular term $t=[\epsilon]$ we have $t={\delta}^i_j{\delta}^j_i$, where $i\not=j$, and for $t=[a]$, where $a\in T$, we have $t=\delta^j_i[a]^i_j$.
For the elementary  singular term $t=[a_{1}\ldots a_{n}]$, where $a_{1},\ldots, a_{n}\in T$ and $n>1$, we have
$t=[a_{1}]^i_{i_1},\cdots,[a_{n}]^{i_{n-1}}_{i}$, where
$i,i_1,\ldots,i_{n-1},$ are pairwise distinct.

For non-elementary terms the statement follows from the above. $\Box$
%It is worth noting that the above defined encoding of terms as context is faithful, i.e. congruence of terms corresponds to equivalence of contexts.
\bp\label{injectivity of tm}
If $t,t'$ are well-formed pseudo-terms and $t\equiv t'$ then the contexts $\mathit{ct}(t)$ and $\mathit{ct}(t')$ are equivalent in the sense that,  for any context $\Gamma$, if the sequents $\vdash\mathit{ct}(t),\Gamma$ and $\vdash\mathit{ct}(t'),\Gamma$ are well-defined then they
  are {\bf TTL} derivable from each other.
\ep
{\bf Proof} Induction on the definition of term equality (i.e. well-formed pseudo-term congruence) using that the context ${\bf\overline a}^i_j,{\bf\overline b}^j_k$ is equivalent to
${\bf\overline a}^i_{j'},{\bf\overline b}^{j'}_k$ and ${\bf\overline a}^i_j$ is equivalent to ${\bf\overline a}^i_k,\bot^k_j$ and to
$\bot^i_k,{\bf\overline a}^k_j$. $\Box$

The two preceding propositions imply that the map $\mathit{ct}$ induces a well-defined map from {\it terms} to equivalence classes of contexts.
\bd
For a tensor term $t$ the notation $\overline{\bf t}$ stands for the context $\mathit{ct}(\tau)$, where $\tau$ is any well-formed pseudo-term representing
$t$. Similarly, the notation ${\bf t}$ stands for the context $\overline{\mathit{ct}(\tau)}$.
The {\it sequent form}  of the classical, respectively intuitionistic typing judgement
$t\rhd A$
is the sequent  $\vdash\overline{\bf t}, A$, respectively ${\bf t}\vdash A$.
\ed
The above definition is, strictly speaking, ambiguous, because the sequent $\vdash\overline{\bf t},A$, respectively
${\bf t}\vdash A$,
   depends on the pseudo-term representing $t$. But by Proposition \ref{injectivity of tm}, different representations of $t$ give rise to equivalent sequents.
   %Our main interest is in regular typing judgements. However, for the completeness theorem that will be given below we will need a more general notion.
%   \bd
%   A term $t$ is
%     is {\it quasi-regular} if it can be written as  $t=[a_1]^{i_1}_{j_1}\ldots[a_n]^{i_n}_{j_n}$, $n>0$, where $a_1,\ldots,a_n$ are pairwise distinct terminal symbols and $i_k\not=j_k$ for $k=1,\ldots,n$.
%
%    A typing judgement $t\rhd A$ is quasi-regular if the term $t$ is quasi-regular.
%    A lexicon is quasi-regular if all its elements are quasi-regular typing judgements.
%   \ed
%   Observe that regular terms are quasi-regular.
   \nb\label{no counits}
   If $t$ is a regular lexical term then there is a representative of $\overline{\bf t}$ that does not contain any counit $\bot^i_j$. $\Box$
   \nbe
%%%%%%%%%%%%%

\subsection{Derivations of typing judgements}\label{typing judgements derivations}
Below we restrict our attention to classical {\bf TTL}.

 Assume that we are given a lexicon $\mathit{Lex}$. We will use the notation ${\bf Lex}$,  ${\bf Lex}_\sim$
${\bf Lex}_{cl}$
for the sets of sequent forms of elements of respectively $\mathit{Lex}$, $\mathit{Lex}_{\sim}$ and $\mathit{Lex}_{cl}$. Note that the sequent forms of typing axioms $\delta^i_j\rhd{\bf 1}^i_j$ are  {\bf TTL} axioms $\vdash\bot_i^j,{\bf 1}_j^i$.  We denote their set as {\bf TpAx}.
\bp\label{cut-elimination for canonical lexicon}
  Let $\mathit{Lex}$ be a  lexicon and assume that the sequent $\vdash\Gamma$ is {\bf TTL}   derivable from elements of  ${\bf Lex}$. Then there exists a derivation $\pi$ of $\vdash\Gamma$ from  elements of ${\bf Lex}$ such that for any application of the Cut rule in $\pi$ one of the two premises
 is an element $\vdash\overline{\bf t}, A$  of ${\bf Lex}$ and   the cut-formula is $A$.
  \ep
  {\bf Proof}
  This follows directly from  Lemma \ref{cut-elimination with axioms}. $\Box$
  \bl\label{deduction theorem}
  Let $\mathit{Lex}$ be a regular lexicon. A sequent  $\Sigma$ of the form $\vdash\overline{\bf t},\Theta$
  (respectively ${\bf t},\Phi\vdash F$ in the intuitionistic case), where $t$ is a tensor term and $\Theta$ is a context in the base language
  (respectively $\Phi$ is a context and $F$ a formula in the base language), is {\bf TTL} (respectively {\bf ITTL}) derivable from elements of ${\bf Lex}$ iff  there exists elements $\vdash\overline{\bf t_{(1)}},A_{(1)},\ldots,\vdash\overline{\bf t}_{(n)},A_{(n)}$
  (respectively ${\bf t_{(1)}}\vdash A_{(1)},\ldots,{\bf t}_{(n)}\vdash A_{(n)}$) of ${\bf Lex}_{cl}$
  such that
  $t=t_{(1)}\cdots t_{(n)}$
  and the sequent
  \be\label{sequent in deduction theorem}
  \vdash\overline{A_{(1)}},\ldots,\overline{A_{(n)}},\Theta
  \ee
  (respectively ${A_{(1)}},\ldots,{A_{(n)}},\Phi\vdash F$)
  is  {\bf TTL} (respectively {\bf ITTL}) derivable.
  \el
  {\bf Proof} For the ``if'' direction use the Cut rule, for the ``only if'' direction use induction on
  a derivation $\pi$ of $\Sigma$ from Proposition \ref{cut-elimination for canonical lexicon}.

  Base of induction: $\Sigma$ is a {\bf TTL} axiom  or an element of ${\bf Lex}$.

  Assume that $\Sigma$ is a {\bf TTL} axiom $\vdash \overline F,F$, where $F\in\mathit{At}$. There are two possibilities. The first one is that $t=1$, so that $\overline{\bf t}=\emptyset$, and $\Theta=\overline F,F$. The statement of the lemma holds with $n=0$. The second possibility is that $\Sigma\in{\bf TpAx}$ is of the form $\vdash\bot^j_i,{\bf 1}^i_j$, the term $t=\delta^i_j$, and the context $\Theta={\bf 1}^i_j$. Then we put $n=1$, $t_{(1)}=t$, and $A_{(1)}={\bf 1}^i_j$. The statement of the lemma obviously holds.

  Assume that the sequent $\Sigma$ is an element of ${\bf Lex}$. This means that $\Sigma$ has the form $\vdash\overline{{\bf s}},A$, where the regular lexical typing judgement $s\rhd A$ is  an element of $\mathit{Lex}$. Since the term $s$ is lexical and regular, by Proposition \ref{no counits}, we may assume that the context $\overline{\bf s}$ does not contain counits. Again,  there are two possibilities. The first one is that $\overline{\bf t}=\overline{{\bf s}}$, $\Theta={A}$, then we put $n=1$, $t_{(1)}=s=t$, $A_{(1)}=A$ and  the statement is obvious. The second possibility is that $\Theta=\emptyset$ and $\overline{\bf t}=\overline{{\bf s}},A$. This may be the case only if $A=\bot_i^j$. Then we have that  $t=s\delta^j_i$. Put $n=2$. $t_{(1)}=s$, $t_{2}=\delta^i_j$, $A_{(1)}=A=\bot^j_i$, $A_{(2)}={\bf 1}^i_j$. The  sequent in (\ref{sequent in deduction theorem}) is nothing but $\vdash\bot^i_j,{\bf 1}^j_i$ and is derivable in {\bf TTL}. The statement holds.

  Now consider the induction steps.

  Assume that $\pi$
   was obtained from derivations $\pi_1,\pi_2$ by the $(\otimes)$ rule. This means that
   that $\Sigma$ has the form $\vdash{\bf \overline r,\overline s},\Gamma,\Delta,F\otimes G$, so that the context
    $\Theta=\Gamma,\Delta,F\otimes G$, the term $t=rs$,
    while $\pi_1$, $\pi_2$ have respectively the conclusions
     $\vdash{\bf \overline r},\Gamma,F$ and $\vdash{\bf \overline s},\Delta,G$.

     By the induction hypothesis there are elements
     $$r_{(1)}\rhd F_{(1)},\ldots,r_{(k)}\rhd F_{(k)},\quad s_{(1)}\rhd G_{(1)},\ldots,s_{(m)}\rhd G_{(m)}$$ of $\mathit{Lex}_{cl}$ such that
     $r=r_{(1)}\ldots r_{(k)}$, $s=s_{(1)}\ldots s_{(m)}$ and the sequents
\be\label{pi_1}
     \vdash\overline {F_{(1)}},\ldots,\overline {F_{(k)}},\Gamma,F,
     \ee
     \be\label{pi_2}
     \vdash\overline {G_{(1)}},\ldots,\overline {G_{(m)}},\Delta,G
      \ee
      are {\bf TTL} derivable.  What we want is to write $t=r_{(1)}\ldots r_{(k)}s_{(1)}\ldots s_{(m)}$, but there is a possibility that the pseudoterm in the righthand side is not well-formed, i.e. there are forbidden index repetitions. Concretely, it may be that some index $i$ occurs twice both in
      the pseudoterm $r_{(1)}\ldots r_{(k)}$ and in the pseudoterm $s_{(1)}\ldots s_{(m)}$.

      Let $I$ be the set of all indices bound in the pseudoterm  $r_{(1)}\ldots r_{(k)}$, and $J$ be the set of all indices bound in $s_{(1)}\ldots s_{(m)}$. Pick fresh disjoint sets $I'$, $J'$ of indices in bijection with $I$, $J$ respectively. Let ${r_{(\alpha)}}'$,  be the term obtained from $r_{(\alpha)}$ by replacing every occurrence of an index $i\in I$ with its image $i'\in I'$, $\alpha=1,\ldots,k$,
     and ${s_{(\alpha)}}'$,  be the term obtained from $s_{(\alpha)}$ by replacing every occurrence of an index $j\in J$ with its image $j'\in J'$, $\alpha=1,\ldots,m$. Similarly,
     let ${F_{(\alpha)}}'$,  be the formula obtained from $F_{(\alpha)}$ by replacing every free occurrence of an index $i\in I$ with its image $i'\in I'$, $\alpha=1,\ldots,k$,
     and ${G_{(\alpha)}}'$,  be the formula obtained from $G_{(\alpha)}$ by replacing every free occurrence of an index $j\in J$ with its image $j'\in J'$, $\alpha=1,\ldots,m$. Then the typing judgements
      $$r_{(1)}'\rhd F_{(1)}',\ldots,r_{(k)}'\rhd F_{(k)}',\quad s_{(1)}'\rhd G_{(1)}',\ldots,s_{(m)}'\rhd G_{(m)}'$$ are elements of $\mathit{Lex}_{cl}$, because $\mathit{Lex}_{cl}$ is closed under similarity. Moreover the term $r_{(1)}'\ldots r_{(k)}'s_{(1)}'\ldots s_{(m)}'$ is well-formed and equal to $rs=t$. Finally, the sequents $$\vdash\overline {F_{(1)}'},\ldots,\overline {F_{(k)}'},\Gamma,F,\quad
     \vdash\overline {G_{(1)}'},\ldots,\overline {G_{(m)}'},\Delta,G$$ are similar to (\ref{pi_1}), (\ref{pi_2}) respectively, hence they are {\bf TTL} derivable, because derivability in {\bf TTL} is closed under similarity. The sequent
     $$\vdash\overline {F_{(1)}'},\ldots,\overline {F_{(k)}'},\overline {G_{(1)}'},\ldots,\overline {G_{(m)}'}\Gamma,
     \Delta,F\otimes G$$ is {\bf TTL} derivable from the above by the $(\otimes)$ rule, and the statement follows.

  Assume that $\pi$
   was obtained from a derivation $\pi'$ by the $(\bot_\leftarrow)$ rule. This means
   that $\Sigma$ has the form $\vdash{\bf s},F_{[\nicefrac{i}{j}}],\Gamma$ while $\pi'$, has the conclusion
     $\vdash{\bf s},\bot^j_i, F,\Gamma$.
     There are two possibilities how the context in $\Sigma$ is partitioned. The first one is that $\Theta=F,\Gamma$ and $t=s$. Then the statement easily follows from the induction hypothesis applied to $\pi'$. The second possibility may occur when
     $F=\bot^k_j$. In this case it may be that
     $\Theta=\Gamma$ and $\overline t=\overline s,F_{\nicefrac{i}{j}}$, i.e. $\overline t=\overline s,\bot^k_i$, so that $t=s\delta^i_k$.
     In this case the conclusion of $\pi'$ is $\vdash{\bf s},\bot^j_i,\bot^k_j,\Gamma$. But the latter sequent can be read as $\vdash\overline{\bf t},\Gamma$, because $s\delta^i_j\delta^j_k=s\delta^i_k=t$. Again, the statement follows from the induction hypothesis.

     Other rules are treated similarly.
   $\Box$
\bc
Given a classical, respectively  intuitionistic, tensor grammar $G=(\mathit{Lex},S)$, a typing judgement $t\rhd F$ is derivable in $G$ iff the sequent $\vdash\overline{\bf t},F$, respectively ${\bf t}\vdash F$, is {\bf TTL}, respectively {\bf ITTL}, derivable from the set of non-logical axioms
$${\bf Lex}=\{\vdash\overline{\bf t}, A|~(t\rhd A)\in\mathit{Lex}\},$$
respectively
$${\bf Lex}=\{{\bf t}\vdash A|~(t\rhd A)\in\mathit{Lex}\}.\quad \Box$$
\ec
%We consider the classical case.
%
%If  $t\rhd_GF$ then the statement follows from the definition of derivability in a grammar (keeping in mind that  the typing axiom $\delta^i_j\rhd{1}^i_j$ corresponds to the {\bf TTL} derivable sequent $\vdash\bot^j_i,{1}^i_j$).
%Assume that  $\vdash\overline{\bf t},F$ is derivable from elements of {\bf Lex}.
%
% %Since derivability
%% in {\bf TTL} or in a tensor grammar
%% is closed under $\alpha$-equivalence, we may assume that  ${\bf Lex}$ and $\mathit{Lex}$ are closed under $\alpha$-equivalence of sequents, respectively typing judgements.
% By Lemma \ref{deduction theorem}, we have a list
% $\vdash\overline{\bf t_{(1)}},A_{(1)},\ldots,\vdash\overline{\bf t_{(n)}},A_{(n)}$ of elements of  ${\bf Lex}$ and indices $i_1,j_1,\ldots,i_k,j_k$ such that
%  $t=\delta^{i_1}_{j_1}\cdots\delta^{i_k}_{j_k}t_{(1)}\cdots t_{(n)}$
%  and  sequent (\ref{sequent in deduction theorem})
%is derivable in {\bf TTL}. Put $t_{(n+\alpha)}=\delta^{i_\alpha}_{j_{\alpha}}$, $A_{(n+\alpha)}={1}^{i_\alpha}_{j_{\alpha}}$, $\alpha=1,\ldots,k$. Then for every $\alpha\leq k$ the typing judgement $t_{(n+\alpha)}\rhd A_{(n+\alpha)}$ is a typing axiom. We have the list
% ${t_{(1)}}\rhd A_{(1)},\ldots,{t_{(n+k)}}\rhd A_{(n+k)}$ of elements of  $\mathit{Lex\cup TpAx}$ such that $t=t_{(1)}\ldots t_{n+k}$ and the sequent
% $\vdash\overline{A_{(1)}},\ldots,\overline{A_{(n+k)}},A$ is derivable in {\bf TTL}. Which means exactly that $t\rhd_GA$.
%
%The intuitionistic version of the statement follows from conservativity of the intuitionistic fragment of {\bf TTL}. $\Box$

Thus we obtain a convenient representation of tensor grammars as systems of non-logical axioms in the extended language.
 When using such a representation  it is natural to add to the system the following
{\it lexical substitution} rule :
$$\cfrac{\vdash A,\Gamma\quad (t\rhd A)\in\mathit{Lex}}{\vdash \overline{\bf t},A}~({\rm{Lex}})$$
in the classical format or
$$\cfrac{A,\Gamma\vdash B\quad (t\rhd A)\in\mathit{Lex}}{{\bf t},A\vdash B}~({\rm{Lex}})$$
in the intuitionistic format. The rule is just a special instance of the Cut rule where one of the premises is an axiom from {\bf Lex}.

\subsection{Example}
 \begin{figure}%[htb]
\centering
%%%%%%%%%%%%%%%%%%%%%%%%%%%%%%%%%
\subfloat[{Axioms}
\label{Axioms}
]
{
$\begin{array}{c}
[{\rm{Mary}}]^i_j\rhd np^i_j,\quad
[{\rm{John}}]^i_j\rhd np^i_j,\quad
[{\rm{loves}}]^i_j\rhd np_y^j\multimap np_i^x\multimap s^x_y,
\\
\mbox{[{\rm{madly}}]}^y_k\rhd
(np^x_i\multimap s^x_y)
\multimap np_i^z\multimap s^z_k,
\quad%\\
\mbox{[{\rm{who}}]}^t_z\rhd
{1}^j_y\otimes(np^j_y\multimap s^z_k)\multimap np_t^u\multimap np^u_k
\end{array}$
}\\
%%%%%%%%%%%%%%%%%%%%%%%%%%%
\subfloat[{Derivation}
\label{ND grammar derivation}
]
{
%$\mbox{ }$\\
 %  \begin{prooftree}
 \begin{tabular}{c}
    \def\ScoreOverhang{.1pt}
    \AxiomC{$
           np^j_y \vdash np^j_y
                $}
     \AxiomC{$ np^x_i \vdash np^x_i$}
     \AxiomC{$ s^x_y \vdash s^x_y$}
    \RightLabel{$(\multimap{\rm {L}})$}
    \insertBetweenHyps{\hskip 5pt}
    \BinaryInfC{$
        np^x_i, np^x_i\multimap s^x_y \vdash s^x_y
        $}
     \RightLabel{$(\multimap{\rm {L}})$}
    \insertBetweenHyps{\hskip 5pt}
    \BinaryInfC{$
        np^x_i,  np^j_y ,np^j_y\multimap np^x_i\multimap s^x_y \vdash s^x_y
        $}
    \RightLabel{$({\rm {Lex}})$}
    \UnaryInfC{$
        np^x_i,   np^j_y ,{\bf loves}^i_j \vdash s^x_y
        $}
    \RightLabel{$(\multimap{\rm {R}})$}
    \UnaryInfC{$
         np^j_y,  {\bf loves}^i_j \vdash np^x_i \multimap s^x_y
        $}
    \AxiomC{$ {\bf John}^z_i \vdash np^z_i$}
     \AxiomC{$ s^z_k \vdash s^z_k$}
    \RightLabel{$(\multimap{\rm {L}})$}
    \insertBetweenHyps{\hskip 5pt}
    \BinaryInfC{$
        {\bf John}^z_i, np^z_i\multimap s^z_k \vdash s^z_k
        $}
    \RightLabel{$(\multimap{\rm {L}})$}
    \insertBetweenHyps{\hskip 5pt}
    \BinaryInfC{$
    {\bf John}^z_i, np^j_y,  {\bf loves}^i_j,( np^x_i \multimap s^x_y)\multimap np^z_i\multimap s^z_k \vdash s^z_k
    $}
    \RightLabel{$({\rm {Lex}})$}
    \UnaryInfC{$
    {\bf John}^z_i, np^j_y,  {\bf loves}^i_j,{\bf madly}^y_k \vdash s^z_k
    $}
     \RightLabel{$(\multimap{\rm {R}})$}
     \UnaryInfC{$
         {\bf John}^z_i, {\bf loves}^i_j,{\bf madly}^y_k \vdash  np^j_y\multimap s^z_k
        $}
         \DisplayProof
\\
%$[\mbox{John}]^z_i[\mbox{loves}]^i_j[\mbox{madly}]^y_k=[\mbox{John loves}]^z_j[\mbox{madly}]^y_k$
\\
    \def\ScoreOverhang{.1pt}
    \AxiomC{$ {1}_y^j \vdash {1}_y^j$}
    \AxiomC{$
            {\bf John}^z_i, {\bf loves}^i_j,{\bf madly}^y_k \vdash  np^j_y\multimap s^z_k
            $}
    \RightLabel{$(\otimes{\rm {R}})$}
    \insertBetweenHyps{\hskip 5pt}
    \BinaryInfC{$
    {\bf John}^z_i, {\bf loves}^i_j,{\bf madly}^y_k,{1}_y^j \vdash{1}_y^j\otimes( np^j_y\multimap s^z_k)
            $}
     \RightLabel{$({\rm {L{1}^\leftarrow}})$}
     \UnaryInfC{$
        {\bf John}^z_i, {\bf loves}^i_j,{\bf madly}^j_k \vdash{1}_y^j\otimes( np^j_y\multimap s^z_k)
        $}
         \DisplayProof
\\
\\
    \def\ScoreOverhang{.1pt}
    \AxiomC{$
        {\bf John}^z_i, {\bf loves}^i_j,{\bf madly}^j_k \vdash{1}_y^j\otimes( np^j_y\multimap s^z_k)
        $}
    \AxiomC{$ {\bf Mary}^u_t \vdash np^u_t$}
     \AxiomC{$ np^u_k \vdash np^u_k$}
    \RightLabel{$(\multimap{\rm {L}})$}
    \insertBetweenHyps{\hskip 5pt}
    \BinaryInfC{$
        {\bf Mary}^u_t, np^u_t\multimap np^u_k \vdash np^u_k
        $}
    \RightLabel{$(\multimap{\rm {L}})$}
    \insertBetweenHyps{\hskip 5pt}
    \BinaryInfC{$
    {\bf Mary}^u_t,{\bf John}^z_i, {\bf loves}^i_j,{\bf madly}^j_k,{1}_y^j\otimes( np^j_y\multimap s^z_k)
    \multimap np^u_t\multimap np^u_k \vdash np^u_k
    $}
    \RightLabel{$({\rm {Lex}})$}
    \UnaryInfC{$
    {\bf Mary}^u_t,{\bf John}^z_i, {\bf loves}^i_j,{\bf madly}^j_k,{\bf who}^t_z \vdash np^u_k
    $}
   \RightLabel{$(=)$}
    \UnaryInfC{$
    {\bf Mary}^u_t,{\bf who}^t_z,{\bf John}^z_i, {\bf loves}^i_j,{\bf madly}^j_k \vdash np^u_k
    $}
  \DisplayProof
\end{tabular}
}
\caption{Tensor grammar example}
\label{ND grammar example}
\end{figure}
A linguistically motivated example of tensor grammar is given in Figure \ref{ND grammar example}.  Figure \ref{Axioms} shows the lexicon, and in Figure \ref{ND grammar derivation} we derive the noun phrase ``Mary who John loves madly'' that we discussed in Section \ref{first label}, using  the representation in the extended language. (For readability, the derivation is broken into three parts.)

\section{Completeness of tensor term semantics}
In this section we establish completeness of {\bf TTL} and, consequently, of {\bf ITTL} in a balanced language for the tensor term semantics. We will discuss the classical system, the intuitionistic case will follow by conservativity.
We assume that we are given a  classical tensor language generated by a finite set $\mathit{Lit}_+$ of balanced positive literals. We also assume that we have an {\it infinite} alphabet $T$ of terminal symbols. We will extensively use the {\it extended language} as discussed in Sections \ref{extended language}, \ref{typing judgements derivations}.

\subsection{Canonical lexicon}
\bd\label{canonical lexicon}
A (classical) lexicon $\mathit{Lex}$  is   {\it canonical} if the following hold:
\begin{enumerate}[(i)]
  \item for any formula $A$  of the base language there are infinitely many elements of  the form $t\rhd A$ in  $\mathit{Lex}$,
    \item if $t\rhd A$, $s\rhd A'$ are distinct elements of $\mathit{Lex}$, then $t$ and $s$ have no common terminal symbols,
        \item 
              For any typing judgement $t\rhd A$
        in $\mathit{Lex}$ the term $t$ can be written as  $t=[a_1]^{i_1}_{j_1}\ldots[a_n]^{i_n}_{j_n}$, $n>0$, where $a_1,\ldots,a_n$ are pairwise distinct terminal symbols and $i_k\not=j_k$ for $k=1,\ldots,n$.
\end{enumerate}
\ed
The elements of a canonical lexicon will serve as ``axioms'' describing the model.
\bp
Canonical lexicons exist.
\ep
{\bf Proof}
A canonical lexicon  can be constructed in countably many steps as follows. Choose an enumeration of the set $\mathit{Fm}(\mathit{Lit}_+)\times\mathbb{N}$.
Let $\mathit{Lex}_0=\emptyset$.
On the $i$th step consider the $i$th element $(A,n)$. Let  $FI(A)=(\{i_1,\ldots,i_k\},\{j_1,\ldots,j_k\})$.  If $k>0$ , choose pairwise distinct terminal symbols $a_1,\ldots,a_k\in T$
 not  occurring  in $\mathit{Lex}_{i-1}$ and define $\mathit{Lex}_i$ as $\mathit{Lex}_{i-1}$ plus the  typing judgement
$[a_1]^{i_1}_{j_1}\ldots[a_k]^{i_k}_{j_k}\rhd A$. If $k=0$, choose $a,b\in T$ not  occurring in $\mathit{Lex}_{i-1}$ and define
$\mathit{Lex}_i$ as $\mathit{Lex}_{i-1}$ plus the  typing judgement
 $[ab]\rhd A$. (Keeping in mind that $[ab]$ can be written as the product $[ab]=[a]^i_j[b]^j_i$.)
 Finally, define $\mathit{Lex}=\bigcup\limits_i\mathit{Lex}_i$.
$\Box$

Henceforth we assume that a canonical lexicon $\mathit{Lex}$ is given. Recall that ${\bf Lex}$ is the set of sequent forms of elements of $\mathit{Lex}$.
 \bp\label{derivable sequent in the base lang}
 If $\vdash\Gamma$ is a sequent in the base language {\bf TTL} derivable from elements of $\mathbf{Lex}$ then $\vdash\Gamma$ is derivable in {\bf TTL}.
 \ep
 {\bf Proof} The statement is a trivial case of Lemma \ref{deduction theorem}. $\Box$
 \bp\label{participation of terminal symbols}
 Let $\vdash\overline{\bf t},\Theta$, where $\Theta$ is a context in the base language, be a sequent {\bf TTL} derivable from elements of  $\mathbf{Lex}$,
  and let  $a\in T$.
  Then  $a$ occurs in $t$  iff there is an element $s\rhd A$ of $\mathit{Lex}_{cl}$ such that $a$ occurs in $s$ and $t$ factors as the product $t=sr$ for some term $r$.
 \ep
 {\bf Proof} This is a direct corollary of Lemma \ref{deduction theorem}. $\Box$
  \bp\label{main for orthogonality}
 Let $s\rhd A$ be an element of  $\mathit{Lex}_{cl}$,  and $t$ be a term having no common terminal symbol with $s$. Assume that the sequent
 $\vdash\overline{\bf t},\overline{\bf s},\Theta$, where $\Theta$ is a context in the base language, is {\bf TTL} derivable from elements of $\mathbf{Lex}$. Then the sequent $\Sigma$ of the form $\vdash\overline{\bf t},\overline A,\Theta$ is also {\bf TTL} derivable from elements of $\mathbf{Lex}$.
 \ep
 {\bf Proof}
 If $s\rhd A$ is a typing axiom then $s=\delta^i_j$, $A={\bf 1}^i_j$ and ${\bf s}=A$, so the statement is obvious.
 So we assume that $s\rhd A$ is an element of $\mathit{Lex}_\sim$.

 Use induction on derivation  $\pi$ of the sequent $\vdash\overline{\bf s},\overline{\bf t},\Theta$ from Proposition \ref{cut-elimination for canonical lexicon}.

 The essential step is when $\pi$ is obtained from derivations $\pi_1,\pi_2$ by the $(\otimes)$ rule. This means that $\Theta$ had the form $\Theta=\Gamma,\Delta,F\otimes G$ and $st$ factors as the product $st=pq$, while $\pi_1,\pi_2$ have respectively the conclusions $\vdash\overline{\bf p},\Gamma,F$ and
 $\vdash\overline{\bf q},\Delta,G$.

  By definition of a canonical lexicon, the term $s$ has representation  $$s=[a_{(1)}]^{i_1}_{j_1}\cdots[a_{(n)}]^{i_n}_{j_n}, \quad\overline{\bf s}=\overline{{\bf a}_{(1)}}_{i_1}^{j_1},\ldots,\overline{{\bf a}_{(n)}}_{i_n}^{j_n},$$ where $a_{(1)},\ldots,a_{(n)}\in T$ are pairwise distinct.

 The terminal symbol $a_{(1)}$ must occur in $p$ or $q$. Assume for definiteness that it occurs in $p$. It follows from Proposition \ref{participation of terminal symbols} and property (ii) of Definition \ref{canonical lexicon} of a canonical lexicon then that all $a_{(1)},\ldots,a_{(n)}$ occur in $p$, hence the context $\overline{\bf s}$ is contained in $\overline{\bf p}$, which means $\overline{\bf p}=\overline{\bf s},\overline{\bf r}$ for some term $r$. Observe that $r$ does not share any terminal symbol with $s$. Indeed, otherwise $r$ would contain all $a_{(1)},\ldots,a_{(n)}$, and then the term $ts=srq$ would contain at least two copies of each of $a_{(1)},\ldots,a_{(n)}$, which contradicts the condition that $t$ has no  terminal symbols in common with $s$.

 Thus the induction hypothesis applies to $\pi_1$ and the sequent $\Pi$ of the form $$\vdash\overline{\bf r},\overline A,\Gamma,F$$ is {\bf TTL} derivable from elements of {\bf Lex}.
 On the other hand, since we have  the equality of contexts $\overline{\bf s},\overline{\bf t}=\overline{\bf p},\overline{\bf q}$ and $t$ has not terminal symbol in common with $s$, it follows that the context $\overline{\bf t}$ is the part of the context $\overline{\bf p},\overline{\bf q}$ that is disjoint from $\overline{\bf s}$, i.e. $\overline{\bf t}=\overline{\bf r},\overline{\bf q}$. Then the desired sequent $\Sigma$ had the form $$\vdash\overline{\bf r},\overline{\bf q},A,\Gamma,\Delta,F\otimes G$$ and is derivable from $\Pi$ and the conclusion of $\pi_2$ by the $(\otimes)$ rule.

 Other induction steps are easy. $\Box$

 \subsection{Canonical model}
\bd
Given a canonical lexicon $\mathit{Lex}$, the {\it canonical type system}  is defined by the dualizing type
\be\label{bbot}
\bbot=
\{t|~\vdash\overline{\bf t}\mbox{ is derivable from }\mathbf{Lex}\}
\ee
The {\it canonical  model} of {\bf TTL} is defined by  the  valuation
 \be\label{canonical valuation}
||A||=\{t|~\vdash\overline{\bf t},A\mbox{ is derivable from }\mathbf{Lex}\},
\ee
where
$A=p^{i_1\ldots i_k}_{j_1\ldots j_k}$, $p\in\mathit{Lit}_+$, $v(p)=(k,k)$.
\ed
It is easy to check that valuation (\ref{canonical valuation}) satisfies the conditions for a model of {\bf TTL} (Definition \ref{valuation}).
%\bp\label{gamma in bbot}
%In the setting of the canonical type system we have $\mathit{tm}(\gamma)\in\bbot$ iff $\vdash\gamma$ is derivable from $\mathit{Lex}$.
%\ep
%{\bf Proof} Assume that $\mathit{tm}(\gamma)=t$ and $t\in\bbot$. This means that there exists some $\gamma'$ such that the sequent $\vdash\gamma'$ is derivable from $\mathit{Lex}$ and $\mathit{tm}(\gamma')=t$. Then $\mathit{tm}(\gamma)=\mathit{tm}(\gamma')$ and by Proposition \ref{injectivity of tm} the sequent $\vdash\gamma$ is derivable from $\vdash\gamma'$. In particular, $\vdash\gamma$ is derivable from $\mathit{Lex}$ as well. The other direction is just definition of $\bbot$. $\Box$
\bl
In the setting of the canonical model, equation (\ref{canonical valuation}) holds for any formula $A$ of the base language.
\el
{\bf Proof} Induction on $A$, where we use that any formula of the base language  can be constructed from  atomic ones of the form
$A=p^{i_1\ldots i_k}_{j_1\ldots j_k}$, where $p\in\mathit{Lit}_+$, and $A=\bot^i_j$ by means of the operations $\overline{(.)}$ and $\parr$.

Base: we only need to consider the case $A=\bot^i_j$. If the sequent $\vdash\overline{\bf t},\bot^i_j$ is derivable from $\mathbf{Lex}$ then, by definition (\ref{bbot}), we have that
$t\delta^j_i\in\bbot$, hence $t\in\overline{\{\delta^j_i\}}=\bot^i_j=||\bot^i_j||$. If, on the other hand,  $t\in||\bot^i_j||=\bot^i_j$ then $t\delta^j_i\in\bbot$, which, again by definition (\ref{bbot}), means precisely that the sequent $\vdash\overline{\bf t},\bot^i_j$ is derivable from $\mathbf{Lex}$. The base of induction is proven.

Let $A=\overline{B}$, and assume that the statement holds for $B$. Assume that the sequent $\vdash\overline{\bf t},A$ is derivable from $\mathbf{Lex}$. Pick any
$s\in ||B||$. By the induction hypothesis, we have that the sequent $\vdash\overline{\bf s},B$ is derivable from $\mathbf{Lex}$.
Using the Cut rule we get that the sequent $\vdash\overline{\bf t},\overline{\bf s}$, i.e. $\vdash\overline{\bf t\cdot s}$, is derivable from $\mathbf{Lex}$ as well, hence $ts\in\bbot$. Since $s\in ||B||$ was arbitrary we conclude that  $t\in \overline{\overline{||B||}}=||A||$. Now assume that $t\in ||A||$.
It follows from Definition \ref{canonical lexicon} that there exists an element
$\vdash\overline{\bf s},B$ in $\mathbf{Lex}$ such that $t$ and $s$ have no common literals. By the induction hypothesis $\overline{\bf s}\in||B||$ and it follows  that $ts\in\bbot$. This means that  the sequent $\vdash\overline{\bf t\cdot s}$, i.e. $\vdash\overline{\bf t},\overline{\bf s}$, is derivable from $\mathbf{Lex}$. Then by Proposition \ref{main for orthogonality} the sequent $\vdash\overline{\bf t},\overline B$, i.e. $\vdash \overline{\bf t},A$ is derivable from $\mathbf{Lex}$ as well.

 The step $A= B\parr  C$
 is very similar to the preceding one.
  %Assume that $\vdash t:A$. If  $s\in\overline{||B||}$, $r\in\overline{||C||}$ then, by the induction hypothesis, $\vdash_Gs:\overline B$, $\vdash_Gr:\overline C$ and using the Cut rule we get that $\vdash_Gr\cdot(rs)$. This means precisely that $t\in ||A\parr B||$. Assume that $t\in ||A\parr B||$. Choose atomic terms $s$, $r$ of types $\overline B$, $\overline C$ respectively such that $t$, $s$, $r$ pairwise have no common terminal symbols. By the induction hypothesis $s\in \overline{||B||}$, $r\in \overline{||C||}$,  hence
%$trs\in\bot$ i.e.
%$\vdash_Gtrs$. Again, from  Proposition \ref{auxiliary for completeness} we  get that $\vdash_Gt:A\parr B$.
%
%Other cases are proven similarly using the Cut rule for one direction and Proposition \ref{auxiliary for completeness} for the opposite one.
 $\Box$
%
%Finally, let $A=\bot^i_j$. If $\vdash_G\sigma:\bot^i_j$ then cutting with the $({1})$ axiom we get that $\vdash_G\sigma\delta^i_j$, hence $\sigma\in \overline{\{\delta^i_j\}}=\bot^i_j$. If $\sigma\in \bot^i_j $ then choose an appropriate atomic term $\tau$ of type ${1}^i_j$ and use
%Proposition \ref{auxiliary for completeness} {\it (ii)} as usual.
%$\Box$
\bc
The sequent $\vdash\Gamma$ in the base language is valid in a canonical model if it is derivable in {\bf TTL}.
\ec
{\bf Proof} Let $\Gamma=A_{(1)},\ldots,A_{(n)}$. The sequent $\vdash\Gamma$ is valid in the canonical model iff  $1\in||A_{(1)}\parr\ldots\parr A_{(n)}||$. By the preceding lemma, the latter condition implies that $\vdash A_{(1)}\parr\ldots\parr A_{(n)}$ is {\bf TTL} derivable from the canonical lexicon. The latter sequent is equivalent to $\vdash\Gamma$, hence $\vdash\Gamma$ is derivable from the canonical lexicon. Then, by Proposition \ref{derivable sequent in the base lang} the sequent $\vdash\Gamma$ is derivable in {\bf TTL}.  $\Box$
\bc[Completeness Theorem]
A sequent in a balanced classical, respectively intuitionistic tensor language is derivable in {\bf TTL}, respectively {\bf ITTL}
iff it is valid in any tensor term model. $\Box$
\ec

\section{Conclusion and future work}
We proposed a simple logical system, {\bf TTL}, which combines commutative and non-commutative operations of linear logic, and
is designed specially for defining categorial grammars.
{\bf TTL}  contains cyclic logic and Lambek calculus (allowing empty antecedents) as conservative fragments. It also allows encoding derivations of ordinary (commutative) multiplicative linear logic. The logic is cut-free and decidable, it has intuitionistic and classical versions and is equipped with a simple intuitive semantics, which is sound and complete. We also defined categorial grammars based on {\bf TTL}.

A number of topics are left for future research. Among them are  complexity of the calculus, algorithms for proof-search and grammar parsing, proof-nets. We plan to study representation of known grammatical formalism, such as {\bf ACG} and {\bf HTLCG}, in {\bf TTL}. Also, connections with first order linear logic remain to be understood. We also think that there are some interesting extensions and modifications of the calculus. Can we add additive connectives? Can we add some structure that will allow representing {\it ordered} tuples of strings and obtain a conservative representation of displacement calculus? We think that other research questions might arise as well, because the setting seems to us rather rich.

\bibliographystyle{plain}
\bibliography{TTL_bibliography}

 \end{document}